\documentclass[12pt]{amsart}
\usepackage{amscd,amsmath,amssymb,amsfonts}

\theoremstyle{plain}
\newtheorem{thm}{Theorem}
\newtheorem{lem}[thm]{Lemma}

\newtheorem{prop}[thm]{Proposition}
\newtheorem{remark}[thm]{Remark}
\newtheorem{remarks}[thm]{Remarks}
\newtheorem{defn}[thm]{Definition}

\newtheorem{conj}[thm]{Conjecture}

\numberwithin{thm}{section}
\numberwithin{equation}{section}

\newcommand{\uh}{\underline{\underline{H}}}
\newcommand{\ldr}{\langle D \rangle}

\newcommand{\p}{\partial}

\newcommand{\eq}[2]{\begin{equation}\label{#2}#1 \end{equation}}
\newcommand{\ml}[2]{\begin{multline}\label{#2}#1 \end{multline}}

\newcommand{\inj}{\hookrightarrow}

\newcommand{\rank}{{\rm rank}}

\newcommand{\Spec}{{\rm Spec \,}}

\newcommand{\OK}{\Omega^1_K}

\newcommand{\OCK}{\Omega^1_{C/K}}

\newcommand{\OG}{\Omega^1_G}
\newcommand{\OGK}{\Omega^1_{G/K}}
\newcommand{\GOG}{\Gamma(G, \Omega^1_G)}
\newcommand{\GOGK}{\Gamma(G, \Omega^1_{G/K})}

\newcommand{\oxx}{{\Omega^2_X}}

\newcommand{\oxxd}{{\Omega^2_X\langle D\rangle(\sD')}}
\newcommand{\rt}{{\rm res \ Tr}}
\newcommand{\trace}{{\rm Tr}}
\newcommand{\sym}{\text{Sym}}
\newcommand{\sC}{{\mathcal C}}
\newcommand{\sD}{{\mathcal D}}
\newcommand{\sE}{{\mathcal E}}

\newcommand{\sG}{{\mathcal G}}

\newcommand{\sI}{{\mathcal I}}

\newcommand{\sK}{{\mathcal K}}
\newcommand{\sL}{{\mathcal L}}

\newcommand{\sO}{{\mathcal O}}

\newcommand{\sS}{{\mathcal S}}

\newcommand{\A}{{\mathbb A}}

\newcommand{\C}{{\mathbb C}}
\newcommand{\D}{{\mathbb D}}

\newcommand{\F}{{\mathbb F}}
\newcommand{\G}{{\mathbb G}}
\renewcommand{\H}{{\mathbb H}}

\newcommand{\N}{{\mathbb N}}

\newcommand{\Q}{{\mathbb Q}}

\newcommand{\V}{{\mathbb V}}
\newcommand{\W}{{\mathbb W}}

\newcommand{\Z}{{\mathbb Z}}

\begin{document}

\title[Gau\ss-Manin Determinants]{Gau\ss-Manin Determinants for Rank
$1$ Irregular Connections on Curves} 
\author{Spencer Bloch}
\address{Dept. of Mathematics,
University of Chicago,
Chicago, IL 60637,
USA}
\email{bloch@math.uchicago.edu}

\author{H\'el\`ene Esnault}
\address{Mathematik,
Universit\"at Essen, FB6, Mathematik, 45117 Essen, Germany}
\email{esnault@uni-essen.de}
\date{August 19, 1999, December 9, 2000}
\begin{abstract}
Let $f:U \to \Spec(K)$ be a smooth open curve over a field
$K\supset k$, where $k$ is an algebraically closed field of
characteristic 0. Let
$\nabla : L \to L\otimes \Omega^1_{U/k}$ be a (possibly irregular)
absolutely integrable connection on a line bundle $L$. A formula
is given for the 
determinant of de Rham cohomology with its Gau\ss-Manin connection
$\Big(\det Rf_*(L\otimes\Omega^1_{U/K}),\ \det\nabla_{GM}\Big)$.  The
formula is expressed as a norm from the curve of a cocycle with values
in a complex defining algebraic differential characters \cite{E}, and this
cocycle is shown to exist for connections of arbitrary rank. 
\end{abstract}
\maketitle
\begin{quote}
Thus mathematics may be defined as the subject in which we never know what we
are talking about, nor whether what we are saying is true.

\vspace*{.5cm}
\noindent Bertrand Russell
\end{quote}

\section{Introduction}
Let $f:U \to \Spec(K)$ be a smooth open curve  
over a field $K\supset k$, where $k$ is an algebraically closed field of
characteristic zero. Let
$\nabla : L \to L\otimes \Omega^1_{U/k}$ be a possibly irregular
absolutely integrable (or vertical, see definition \ref{def2.15})
 connection on a line bundle $L$.
The Riemann-Roch problem in this context is
to describe characteristic classes for the relative de Rham cohomology
$Rf_{*}(L\otimes\Omega^*_{U/K})$ as a (virtual) vector space over
$K$ with an integrable connection, in terms of data on $U$. The $0$-th
characteristic class, the Euler characteristic $\dim R^0 - \dim R^1$,
is well-known to be given by 
\eq{2-2g-n-\sum_i \max(0,m_i-1)
}{1.1}
where $g$
is the genus of the complete curve $C$, $n$ is the number of missing
points, and $m_i$ is the order of the polar part of the connection at
the $i$-th missing point. The purpose of this article is to give a
formula for the first characteristic class, which is the determinant
of the Gau\ss-Manin connection on the relative de Rham cohomology of
the line bundle,
\eq{\Big(\det(Rf_*(L\otimes\Omega^*_{U/K})), \det\nabla_{GM}\Big).
}{1.2}
When $U=C$, so the connection has no poles, the formula given in
\cite{BE} is
\eq{\Big(\det(Rf_*(L\otimes\Omega^*_{C/K})), \nabla_{GM}\Big) =
-f_*((L,\nabla)\cdot 
c_1(\OCK)). 
}{1.3}
Concretely, if one has $c_i \in C(K)$ with $\sum c_i$ a $0$-cycle in
the linear series representing 
$\OCK$, then the determinant is given by restricting $L$ with its
connection to each $c_i$ and then tensoring the resulting lines
with connection together.  

When the connection $\nabla$ has at worst regular singular points at
the points in $D:=C-U$ there is an analogous formula using linear
series given by 
divisors of rational sections $s$ of $\OCK(D)$ satisfying the rigidity
condition ${\rm res}_D(s) = 1$. Indeed, these formulas are valid also 
for higher rank connections. One takes the determinant at zeroes and poles of
$s$. 

In the case of irregular singular points, a
similar formula is possible, but the
rigidification taken must depend on the polar part of the connection. 
Let $(\sL, \nabla)$ be an extension of $(L,
\nabla)$ to $C$,  $\sD= \sum_i m_iD_i$ be a divisor with
multiplicities $m_i \ge 1$ supported in $C-U$ such that the relative
connection 
\eq{\nabla_{/K}:\sL \to \sL\otimes\OCK(\sD)
}{1.4}
yields a complex quasiisomorphic to
$j_*L \to j_*L\otimes \Omega^1_{(C-D)/K}$ and has poles
at all points $D_i$. Then 
$\nabla_{/K}$
does not factor through
\eq{\nabla_{/K}:\sL \to \sL\otimes\OCK(\sD-D_i)}{1.5}
for any $i$.
Writing $\sD$ also for
the artinian subscheme of $C$ determined by $\sD$, this implies
that $\nabla_{/K}$ induces a {\it function linear isomorphism}
\eq{\nabla|_\sD : \sL|_\sD \stackrel{\cong}{\to} \sL\otimes
\OCK(\sD)|_\sD}{1.6}
Because these maps are function linear, we may cancel the lines
$\sL|_\sD$ and deduce canonical elements ${\rm triv}_\nabla \in
\OCK(\sD)|_\sD$. We view ${\rm triv}_\nabla$ as a trivialization of
$\OCK(\sD)$ along $\sD$.  It is known (\cite{EV}, Appendix B) 
that the coboundary of
${\rm triv}_\nabla$ in $H^1(C,\omega_{C/K}) \cong K$ is given by the
degree of $\sL$. Our main result is:
\begin{thm} Let notation be as above. Assume all $D_i$ are defined
over $K$ and some $m_i \geq 2$. Because we are only concerned with the
cohomology over $X-\sD$, we can take
$\sL$ of degree 0 so that the element ${\rm triv}_\nabla$ can be lifted
to $ H^0(C,\OCK(\sD))$. Let $s$ be any such lifting, and write $(s)$
for the divisor of $s$ as a section of $\OCK(\sD)$. Note the support
of (s) is disjoint from $\sD$. Then 
\ml{\det(Rf_*(L\otimes\Omega^*_{U/K})), \det(\nabla_{GM}) \cong \\
-f_*(L\cdot (s))+\tau(L) \in \Omega^1_K/d \log K^*.
}{1.7}
Here $\tau(L)$ is a
$2$-torsion term which can be written
$$ \tau(L) = \sum_i \frac{m_i}{2}d\log(g_{i,0})\in
\frac{1}{2}d\log(K^\times)/d\log(K^\times)
$$ 
where the connection $\nabla_{/K} =
(g_{i,0}+g_{i,1}z_i+\ldots)dz_i/z_i^{m_i}$ for a local coordinate $z_i$
at $D_i \in \sD$. 
\end{thm}
Note that $\Omega^1_K/d\log K^\times$ is the group of isomorphism classes
of rank $1$ connections on $\Spec(K)$. Our assumption that points of
$\sD$ are defined over $K$ is made to avoid complications involving
generalized jacobians in \S 2.
 We remark, of course, that part of
our task will be to give a precise definition of the right hand side of
the  formula of the theorem. It will appear as a product followed by a
trace, and this definition does not depend on the particular choice of
$\sL$ above. In particular, this gives a formulation if we don't assume
that 
$\sL$ is of degree 0, and also if we don't assume that
$m_i \geq 2$ for at
least one $i$, that is if $\nabla$ has regular singular points (see theorem
\ref{thm4.5}). The precise general formulation of our theorem is in
\ref{thm4.7}.  In the case
that $(s)$ is a sum of $K$-points $c_i$, one may simply take the tensor
product of the lines with connection $L|_{c_i}$. 
The right hand side of the formula depends only on the equivalence
class of $(s)$ in a generalized Picard (or divisor class) group of
line bundles with trivializations along $\sD$. Thus, by analogy with
\eqref{1.3}, it is natural to write formula \eqref{1.7} in the form 
\ml{\big((\det(Rf_*(L\otimes\Omega^*_{U/K})), \det(\nabla_{GM})\Big)
\cong \\ 
f_*\Big(L\cdot c_1(\OCK(\sD),{\rm triv}_\nabla))\Big)^{-1} + \tau(L).
}{1.8}
The classical Riemann-Roch pattern begins to break down in that the
characteristic class $c_1(\OCK(\sD),{\rm triv}_\nabla)$ depends on more
than just the geometry of $f:U \to \Spec(K)$.
This reflects the fact that the de Rham cohomology of an irregular
connection depends on more than topology.

There is an analogy here with the case of $\ell$-adic sheaves. If $\sE$ is an
unramified $\ell$-adic sheaf on a complete curve $C$ over a
finite field $\F_q$, 
then the global epsilon factor is given by
$$(-F_q|\det(H^*_{\text{\'et}}(C_{\overline{\F}},\sE))) = \det(\sE)^{-1}\cdot
c_1(\Omega^1_{C/\F_q}).
$$
The basic result in the ramified $\ell$-adic case (\cite{L}) is that the global
epsilon factor can be written as a product of local terms corresponding to
points on the curve where the sheaf ramifies or where a chosen meromorphic
$1$-form has zeroes or poles.  We suspect formula \eqref{1.8} is
analogous to a classical formula for Gau\ss\ sums 
$$g(c,\psi) = \sum_{a\in (\sO/\frak f)^\times} c(a)\psi(a)
$$
where $\frak f \subset \sO$ is an ideal in the ring of integers in a
local field, $c$ (resp. $\psi$) is a character of $(\sO/\frak
f)^\times$ (resp. $(\sO/\frak f)^+$), and both $c$ and $\psi$ have
conductor $\frak f$. If the residue field of $\sO$ has $q$ elements with
$q$ odd, one finds
$$g(\epsilon,\psi) = \begin{cases}q^nc(x) & \frak f=\frak
m^{2n},\ n\ge 1 \\
q^nc(x)\sigma & \frak f=\frak
m^{2n+1},\ n\ge 1. \end{cases}
$$
In this formula $x\in \sO/\frak f$ is a suitable
point, $\sigma=\zeta\sigma_0$ with
$\zeta^q=1$ and
$\sigma_0^2 =
\Big(\frac{-1}{\F_q}\Big)$. (Here $\sigma_0$ is a quadratic Gau\ss\ sum.)

  Our proof follows the main idea of Deligne
\cite{De}. For computing the $\epsilon$-factor associated to a rank one
Galois representation on a curve, he expresses the determinant
of the cohomology as the cohomology on a symmetric
product of $(C-D)$ and reduces the computation to
the geometry of the generalized jacobian. 
In the geometric situation  one is further able
to express the determinant
Gau\ss-Manin connection as the connection arising by restricting a certain
translation-invariant connection to one specific $K$-point of the generalized
jacobian. The essential point seems to be that the de Rham cohomology of a
connection of the form $d+\omega$ on a trivial bundle is somehow concentrated
at the points where $\omega = 0$. 

In section 4 we reinterpret the Riemann-Roch formula in terms of a pairing
\eqref{4.5}
\begin{gather}
\cup: \H^1(C, \sO^*_C \to \sO^*_{\sD}) \times
\H^1(C, \sO^*_C \to \Omega^1_C\ldr(\sD')) \notag \\
\to \H^2(C, \sK_2 \to \Omega^2_C)
\end{gather}
and a trace map \eqref{4.6}
$${\rm Tr \ }:\H^2(C, \sK_2 \to \Omega^2_C) \to \Omega^1_K/d\log K^*.
$$
In section 5 we give an analogous ``non-commutative''
product formula in the higher rank case which
we conjecture calculates the determinant connection in the generic situation
when the connection defines local isomorphisms $E|_\sD \cong E|_\sD\otimes
\omega_{\sD/K}$ (see \eqref{ass5.3}) and the poles of the
absolute connection behave well (see \eqref{ass5.1}). 
We verify the formula has the appropriate invariance properties.
We also show that there is a more general higher rank product of
which it is a special case.
Finally, in section 6 we give a general formula which calculates the group of
isomorphism classes of irregular, integrable, rank $1$ connections in higher
dimensions on a smooth projective variety. 

 We apologize for not expressing our results in the modern language of
$\sD$-modules, but in fact for the study of Gau{\ss}-Manin determinants
there is little gain in passing from connections to $\sD$-modules. Also,
rigidity for connections means that the Gau{\ss}-Manin 
determinant connection is
determined by its value at the generic point on the base, so we may work
with curves over a function field $\Spec(K)$. 

It is our pleasure to
acknowledge the intellectual debt we owe in this work to P. Deligne. We
are also grateful to the Humboldt foundation for financing which enabled
us to work together. 

\section{Connections and forms on Generalized Jacobians}

Throughout this paper $C$ will be a smooth projective curve over a field
$K$ containing an algebraically closed subfield $k$ of
characteristic 0, and $\sD = \sum m_i c_i$ is a divisor on $C$,
with $c_i \in C(K)$. 
We write $G=J_{\sD}$ for the
generalized Jacobian parametrizing isomorphism classes of degree $0$
line bundles on $C$ with trivialization along $\sD$. Fixing a
$K$-rational point $$c_0 \in (C-D)(K),$$ there is
a cycle map $i: C-D \to J_{\sD}$ associating to a closed point $x\in C$ with
$[K(x):K]=n$ the class of the line bundle $\sO(x-nc_0)$ together with the 
trivialization $b|_\sD\circ(a|_\sD)^{-1}$, where 
$$\sO_C \stackrel{a}{\hookleftarrow} \sO_C(-nc_0) \stackrel{b}{\hookrightarrow}
\sO_C(x-nc_0)$$
are the natural maps.

The aim of this section is to describe invariant line bundles
with connection on $J_{\sD}$, comparing them via the cycle
map $i$ to line bundles with connection on $(C-D)$ with a
certain irregularity behavior along $D$.

When the line bundle in question is the trivial bundle, this amounts to
studying invariant (absolute) differential forms on the generalized
jacobian, so we should start with that. Before doing so, however, it is
necessary to understand global functions on the generalized
jacobian. We write
\eq{G \twoheadrightarrow G_0 \twoheadrightarrow J}{2.1}
where $J$ is the usual Jacobian of $C$, and $G_0$ is a semi-abelian
variety. We have extensions 
\begin{gather}
0 \to T \to G_0 \to J \to 0 \label{2.2}\\
0 \to \V \to G \to G_0 \to 0 \label{2.3}
\end{gather}
Here $\V$ is a vector group (isomorphic to $\Spec(\sym(V^*))$ for some
vector space $V$) and $T$ is a torus, i.e. $T_{\bar{K}}\cong
\G_m^r$. 

\begin{lem}\label{lem2.1}The semi-abelian variety $G_0$ admits a universal
vectorial extension 
\eq{0 \to \W \to \sG \to G_0 \to 0.
}{2.4}
In fact, this extension is given by the pullback to $G_0$ of the
universal vectorial extension over $J$. In particular, $\W =
\Gamma(J^\vee,\Omega^1_{J^\vee/K})\otimes \G_a$. 
\end{lem}
\begin{proof} It will suffice to show the pullback vectorial extension
is universal. Since $\text{Ext}^1(T,\G_a)=(0)=\text{Hom}(T,\G_a)$, any
extension of $G_0$ by a vector group $\W$ is pulled back from a unique
extension of $J$ by $\W$. This extension of $J$ is a pushout from the
universal vectorial extension, so the same holds for the pullbacks to
$G_0$. 
\end{proof}

\begin{lem}\label{lem2.2} Let $\pi:G_0\to J$ be an extension of $J$ by
$T$ as above. There 
exists, possibly after a finite field extension,
a quotient torus $T \twoheadrightarrow S$ and a diagram
\eq{\begin{array}{ccc} \makebox[0cm][c]{$T$} &  \hookrightarrow &
\makebox[0mm]{$G_0$} \\
\ \ \ \makebox[0mm][r]{surj. $\downarrow$ }  & \swarrow a \\
\makebox[0cm][c]{$S$}
\end{array}
}{2.5}
such that 
\eq{H^i(G_0,\sO_{G_0})\cong H^i(J,\sO_J) \otimes_{K} H^0(S,\sO_S)
}{2.6}
\end{lem}
\begin{proof}There is a boundary map
\eq{\p:{\rm Hom}_{\bar{K}}(T,\G_m) \to 
{\rm Ext}^1_{\bar{K}}(J, \G_m)
}{2.7}
Define $N:=\ker(\p) \subset M:=Hom_{\bar{K}}(T,\G_m) $. 
Let $S= \text{Hom}(N,\G_m)$ be the torus with character group
$N$.  
For $m\in M$ let $L(m)$ be the line bundle on $J_{\bar{K}}$
corresponding under the map \eqref{2.7}. As an $\sO_{J_{\bar K}}$-algebra
\eq{\pi_*\sO_{G_{0,\bar K}} \cong \oplus_{m\in M}L(m) \cong
H^0(\sO_S)\otimes \Big(\oplus_{m\in M/N}L(m)\Big)
}{2.8} 
The map $a$ in the diagram \eqref{2.5} comes from the above inclusion
$$H^0(\sO_S)\otimes L(0) \subset \pi_*\sO_{G_0}.$$
For $m\in M/N$, (as is well known, cf. \cite{Mu} III 16), 
$L(m)$ has trivial cohomology in all
degrees unless $m=0$. The proposition follows by taking cohomology of
\eqref{2.8}. 
\end{proof}

\begin{lem}\label{lem2.3} Let notation be as above. Let
\eq{
\begin{CD}0 @>>> H^0(J^\vee,\Omega^1_{J^\vee/K})\otimes \G_a @>>>
\sG @>p>> G_0 @>>> 0\end{CD}
}{2.9} 
be the universal vectorial extension. Then
\eq{
H^0(\sG,\sO_\sG)\cong H^0(G_0, \sO_{G_0})
}{2.10}
\end{lem}
\begin{proof} The $\sO_{G_0}$-algebra $p_*\sO_\sG$ is filtered, with
\eq{
gr_ip_*\sO_\sG = fil_i/fil_{i-1} \cong
\sym^i(H^0(J^\vee,\Omega^1_{J^\vee/K})^*)\otimes \sO_{G_0} 
}{2.11} 
With respect to the exact sequences
\eq{\begin{CD} 0 @>>> fil_{i-1} @>>> fil_i @>>> gr_i @>>> 0  \end{CD}
}{2.12}
it suffices to show the boundary map
\eq{b: \sym^i(H^0(J^\vee,\Omega^1_{J^\vee/K})^*)\otimes H^0(G_0,\sO_{G_0}) \to
H^1(G_0, fil_{i-1}) 
}{2.13}
is injective. Composing on the right with the evident map, it suffices
to show the maps
\ml{\sym^i(H^0(J^\vee,\Omega^1_{J^\vee/K})^*)\otimes H^0(G_0,\sO_{G_0}) \\
\to
H^1(G_0, \sO_{G_0})\otimes 
\sym^{i-1}(H^0(J^\vee,\Omega^1_{J^\vee/K})^*)
}{2.14}
are injective. But
\ml{H^0(J^\vee,\Omega^1_{J^\vee/K})^*\otimes H^0(G_0,\sO_{G_0}) \cong
H^1(J, \sO_J)\otimes H^0(G_0,\sO_{G_0}) \\
\cong H^1(G_0, \sO_{G_0})
}{2.15}
and the map in \eqref{2.14} is the map $x^i\mapsto x\otimes x^{i-1}$,
which is injective.   
\end{proof}

\begin{lem}\label{lem2.4} Let $G=J_\sD$ be a generalized jacobian as
above. Then there exists a commutative affine algebraic group $\G$ over
$K$ and a map $\psi : G \to \G$ such that 
\eq{\psi^* : H^0(\sO_\G) \cong H^0(\sO_G).
}{2.16} 
\end{lem}
\begin{proof}Take $\G = \Spec(H^0(G,\sO_G))$.   
\end{proof}

\begin{lem}\label{lem2.5} Let $A$ be the coordinate ring of a
commutative affine algebraic group $H$ over a field $K$ of
characteristic $0$. Corresponding to the simplicial algebraic group
$BH$, one has a complex
\eq{
\begin{CD} A @>\mu^*-p^*_1 -p^*_2 >> A\otimes_K A @> p_{23}^* -
\mu_{12}\otimes p_3^* + p_1^*\otimes \mu_{23}^* - p_{12}^* >> A\otimes
A\otimes A \end{CD}
}{2.19}
This complex is exact at the middle term. 
\end{lem}
\begin{proof} By proposition 4 on p. 168 of \cite{Se}, 
the cohomology in the middle is a
subgroup of the group of extensions $Ext(H,\G_a)$. (Note, $A
=Map(H,\G_a)$.) By the classification of commutative algebraic groups
in characteristic $0$, this ext group vanishes (cf \cite{Se}, pp. 170-172).  
\end{proof}

We write $\OG$ (resp. $\OGK$) for the sheaf of $1$-forms relative to
$k$ (resp. $K$). 
Now we would like to define invariant bundles, connections, differential forms,
cohomology classes of $\sO_G$. 
\begin{defn}\label{def2.6}\begin{enumerate}
\item A rank one bundle $L \in H^1(G, \sO^*_G)$ is called invariant
if $\mu^*L= p_1^*L \otimes p_2^*L \in H^1(G\times G, \sO^*_{G
\times G})$, where $\mu: G \times G \to G$ is the multiplication
and $p_i: G \times G \to G$ are the projections.
\item A global $1$-form $\eta\in \GOG$ is called
invariant if 
$\eta(0)=0 \in \OK$ and $\mu^*\eta = p_1^*\eta + p_2^* \eta \in
\Gamma(G\times G, \Omega^1_{G\times G})$. 
\item A rank one bundle with a connection $(L, \nabla) \in \H^1(G,
\sO^*_G \to \Omega^1_G)$ is called invariant if $(L,
\nabla)|\{0\}= 0 \in \H^1(\Spec K, \sO^*_{\Spec K} \to
\Omega^1_{\Spec K}) = \Omega^1_K /d\log K^*$ and 
$\mu^*(L, \nabla)= p_1^*(L, \nabla) \otimes p_2^*(L, \nabla)
\in \H^1(G \times G, \sO^*_{G \times G} \to \Omega^1_{G \times
G})$. 
\item A class $s\in H^i(G, \sO_G)$ is called 
invariant if $s|\{0\}=0$ and 
$\mu^* s= p_1^*s + p_2^*s$ in $H^i(G\times G, \sO_{G \times
G})$. 
\end{enumerate}
We denote by $H^1(G, \sO^*_G)^{{\rm inv}}$, 
$\GOG^{{\rm inv}}$, $\H^1(G,
\sO^*_G \to \Omega^1_G)^{{\rm inv}}$, and $H^i(G, \sO_G)^{{\rm inv}}$ 
the corresponding groups of
invariant bundles, forms, connections and classes. 
One defines similarly the
groups of relative invariant forms 
$H^0(G, \Omega^1_{G/K})^{{\rm inv}}$  and relative invariant
connections $\H^1(G, \sO^*_G \to \Omega^1_{G/K})^{{\rm inv}}$
without condition on the restriction to the zero section, and
observe that the natural map
$\Omega^1_G \to \Omega^1_{G/K}$ takes global invariant groups to
relative invariant groups. 
\end{defn}

\begin{remark} In the above definitions, we could have defined a weaker notion
of invariance by allowing constant elements. 
We adopt here the rigidification at the origin, keeping in mind that 
without this condition, the corresponding groups obtained
are a direct sum of the ones obtained with the rigidification
and the value of the group on the zero section. Notice, for example, that with
our definition, nonzero constant functions are not invariant! Of course, for
relative objects, there is no distinction. 
\end{remark}

\begin{lem}\label{lem2.7} Let $G=J_{\sD}$ be a generalized Jacobian as
above. Let 
$\tau \in \GOGK^{{\rm inv}}$ be an invariant relative $1$-form on $G$, and
assume $\tau $ lifts to an absolute global form. Then $\tau$ lifts to
an invariant absolute form on $G$.
\end{lem}
\begin{proof}Let $\eta \in \GOG$ be an absolute lifting. Replacing
$\eta$ with $\eta - \eta(0)$ we may assume $\eta(0)\in \OK$
vanishes. Then
\eq{(\mu^*-p_1^* -p_2^*)(\eta)\in H^0(\sO_{G\times_K G})\otimes\OK
}{2.20}
vanishes in $H^0(\sO_{G\times_K G\times_K G})\otimes\OK$. Let $\psi: G
\to \G$ be as in lemma \ref{lem2.4}, 
so $\psi^*:A := H^0(\sO_{\G}) \cong H^0(\sO_G)$. The
previous lemma implies there exists $\sigma \in H^0(\sO_G)\otimes \OK$
with $(\mu^*-p_1^* -p_2^*)(\sigma) = (\mu^*-p_1^* -p_2^*)(\eta)$. Then
$\eta - \sigma$ is the desired invariant absolute form. 
\end{proof}

We next need to relate connections on the curve with invariant connections
on the generalized Jacobian. Here $G = J_{\sD}$ with $\sD = \sum
m_i c_i$. Also, 
\eq{D := \sum c_i;\quad \sD' := \sD -D.
}{2.21}
We assume $D\neq \emptyset$. 

First, rather briefly, we consider the question of what poles
invariant forms on $G$ have when pulled back to $C-D$ via the jacobian
map $C-D \to G$ (defined once we have a basepoint in $C-D$). For simplicity, we
continue to assume the $c_i$ are defined over $K$. Consider the diagram: 
\eq{\begin{array}{ccccc} G \twoheadrightarrow G_0 
& \twoheadrightarrow & J \\ 
\makebox[1.3cm][l]{$\uparrow i$} & & \makebox[.2cm][l]{$\uparrow i'$} \\
\makebox[0cm][c]{$C-D\quad $} & \stackrel{j}{\hookrightarrow}  & C
\end{array}
}{2.22}
We want to compute the pullbacks $i^*H^0(\OG)^{inv}$ (resp.
$i^*H^0(\OGK)^{inv}$) in $H^0(C,\Omega^1_C(*D))$ (resp. in
$H^0(C,\Omega^1_{C/K}(*D))$) .

Pulling back $G$ via $i'$ we get a torseur  
\eq{{i'}^*G \stackrel{p}{\to} C
}{2.23}
under the group $\sO^*_{\sD}/\G_m = \prod \sO^*_{m_\ell c_\ell}/\G_m$. Fix
$\ell$ and let $R=k[[t_\ell]]\subset M=k((t_\ell))$ where $t_\ell$ is a formal
parameter at 
$c_\ell$.  Fix a splitting of the torseur
\eq{G_R \cong \sO^*_{\sD}/\G_m \times \Spec(R).
}{2.24}
Let $c_0\in (C-D)(K)$ be the base point used to define the
map $i:C-D \to G$. Let $\sL$ on $C\times C$ be a line bundle with
$\sL|_{\{c\}\times C}\cong \sO(c-c_0)$. Note one has 
\eq{\sO_{C\times C}\leftarrow p_2^*\sO(-c_0) \to \sL
}{2.25}
and these maps are isomorphisms on 
\eq{\Spec(M)\times \Spec(\sO_{C,D}) \subset
\Spec(M)\times C
}{2.26} 

Corresponding to $\sL|_{\Spec(M)\times C}$ and the
above trivialization, one gets a map $u: \Spec(M) \to G$. With respect
to the above splitting, we view $u$ as an element 
\eq{\prod_i (u_{i0}+u_{i1}t_i+\ldots+u_{i,m_i-1}t_i^{m_i-1})\in
\prod_i (\sO_{m_ic_i}\otimes_K M)^* \mod M^*
}{2.27}
As described in \cite{Se} VII 4, 21, the local shape of $u$
around $c_\ell$ is given by taking the rational function 
$(s_\ell-t_\ell)^{-1} $, where the local coordinates around
$(c_\ell, c_\ell)$ in $C \times C$ is $(s_\ell, t_\ell)$, and
considering it as a unit in $\sO_{m_\ell c_\ell}\otimes M
\cong M[t_\ell]/<t_\ell^{m_\ell}>$. (We change notation so $s_\ell$ is the
local parameter in $R\subset M$.) Since
$u$ is well defined and non-vanishing in $c_i$ for $i\neq \ell$, we have 
\eqref{2.26} that 
\eq{u_{ij}\in R, u_{i0}\in R^*\text{ if } i\neq \ell,\
\text{ord} (u_{\ell
0})=1}{2.28}

The pullbacks to $\Spec(M)$ of the invariant relative differential
forms on $G$ are given by the pullback of invariant relative forms on
$J$ together with the coefficients of powers of the $T_i \mod
T_i^{m_i} $ in the formal
expression 
\begin{multline}\label{2.29}
\sum_i
(u_{i0}+u_{i1}T_i+\ldots+u_{i,m_i-1}T_i^{m_i-1})^{-1} \times \\
(du_{i0}+ du_{i1}T_i+
\ldots+
du_{i,m_i-1}T_i^{m_i-1})  \\
= \sum_i \tau_{i0}+\ldots + \tau_{i,m_i-1}T_i^{m_i-1}.
 \end{multline}
Then \eqref{2.28} implies that
\begin{gather}\label{2.30}
 \tau_{ij} \in \Omega^1_{R} \ {\rm for \ } i \neq \ell \notag \\
\tau_{\ell j} \in \Omega^1_R\ldr( jD_\ell)-\Omega^1_R\ldr( (j-1)D_\ell). 
\end{gather}
Here we denote by 
$\Omega^1_C\ldr$ the sheaf of absolute differential forms of degree
1 with logarithmic poles along $D$. (See formula \ref{pbtau} for a more
precise computation).

These are not all the
absolutely invariant forms, however. One also has forms pulled back
from $J$, but these are regular along $D$. Finally, from lemma
\ref{lem2.7} one has an exact sequence
\ml{0 \to H^0(\sO_G)^{{\rm inv}} \otimes \OK \to
H^0(G,\OG)^{{\rm inv}}
\to H^0(G,\OGK)^{{\rm inv}} \\
\to H^1(\sO_G)^{{\rm inv}} \otimes \OK
}{2.31}
It shows one must consider invariant forms in
$H^0(\sO_G)^{{\rm inv}}\otimes \OK$. We will see in the proof of
proposition \ref{prop2.12} below that these map to $\OG(\sD')$. In sum,
the above discussion shows that the maps in the following proposition
are defined.
\begin{prop}Pullback gives isomorphisms
\begin{gather}\label{2.32}H^0(\OGK)^{inv}
\stackrel{\cong}{\longrightarrow}  H^0(C, \Omega^1_{C/K}(\sD)) \\ 
\label{2.33} H^0(\OG)^{inv} \stackrel{\cong}{\longrightarrow}  H^0(C,
\Omega^1_C\ldr(\sD')) 
\end{gather}
\end{prop}
\begin{proof} Pullback on invariant relative forms is injective,
because $G$ is 
generated by the image of $C-D$. It follows by dimension count that
the first arrow \eqref{2.32}
above is an isomorphism. For the absolute forms we may
consider the diagram
\begin{tiny}
\begin{equation}\label{2.34}\begin{array}{ccccccccc}
0 & \to & H^0(\sO_G)^{{\rm inv}}\otimes\OK & \to & H^0(\OG)^{{\rm inv}} &
\makebox[0cm][r]{$\to$} & H^0(\OGK)^{{\rm inv}} & \to & \quad\quad\ \
\makebox[1.3cm][r]{$H^1(\sO_G)^{{\rm inv}}\otimes
\Omega^1_K$} \\ 
&& \downarrow \cong && \downarrow && \downarrow \cong &&
\downarrow
\cong \\
0 & \to & \makebox[1cm][c]{$H^0(\sO_C(\sD'))\otimes\OK$} & \to &
\makebox[2.2cm][c]{$\ \ H^0(\Omega^1_C\ldr(\sD'))$}   & \makebox[0cm][r]{$\to$} &
H^0(\Omega^1_{C/K}(\sD)) &\ \makebox[0cm][r]{$\to$} & \quad\quad\ \ 
\makebox[1.3cm][r]{\ \ $H^1(\sO_C(\sD'))\otimes\OK$} 
\end{array}
\end{equation}
\end{tiny}

The left and right hand vertical arrows are shown to be isomorphisms
in the proof of proposition \ref{prop2.12}. Hence the isomorphism on invariant
relative forms implies the isomorphism \eqref{2.33} 
on invariant absolute forms.
\end{proof}

We now consider invariant connections on line bundles on $G$.

\begin{lem}\label{lem2.8}Assume the toric subquotient $T$ of $G$ has
trivial Picard group (e.g. $T$ split). Then the map $\H^1(G, \sO_G^*
\to \OGK)^{{\rm inv}} \to H^1(\sO_G^*)^{{\rm inv}}$ is surjective. 
\end{lem}
\begin{proof}This follows because
\eq{\H^1(J, \sO_J^* \to \Omega^1_{J/K})^{{\rm inv}} \twoheadrightarrow
H^1(\sO_J^*)^{{\rm inv}}\twoheadrightarrow H^1(\sO_G^*)^{{\rm inv}}. 
}{2.35}
The second arrow is surjective because we have a diagram
\eq{\begin{CD}0 @>>> H^1(\sO_J^*)^{{\rm inv}} @>>> H^1(\sO_J^*) @>>>
H^2_{DR}(J/K) \\
@. @VVV @VV surj. V @VV inj. V \\
0 @>>> H^1(\sO_G^*)^{{\rm inv}} @>>a > H^1(\sO_G^*) @>>b >
H^2_{DR}(G/K) \\
\end{CD}
}{2.36}
The bottom row is not a priori exact, but $b\circ a=0$ (because
$(N\delta)^*$ acts by $N^2$ on $H^2_{DR}(G/K)$.) The middle
vertical arrow is onto e.g. because the Picard group of the generic
fibre of $G \to J$ is zero. Indeed, $G\to J$ is rationally split, and
the kernel has trivial Picard group by hypothesis. (Since the function
field of the generic fibre equals the function field of $G$, any
divisor on $G$ can be moved by rational equivalence to avoid the
generic fibre, i.e. to be a pullback from the base.)

Finally, the right hand vertical
arrow is injective because, after making a base change
$K \subset \C$,  one can think of $G$ and $J$ as quotients
of vector spaces by lattices, and the map on lattices is surjective
$\otimes \Q$. 
\end{proof}  
\begin{lem}\label{lem2.10} Let $a\in \H^1(G, \sO^*_G \to \OGK)^{{\rm inv}}$
be an invariant connection on a line bundle on the generalized jacobian
$G$. Suppose $a$ lifts to an absolute connection $b'\in \H^1(G,
\sO^*_G \to \OG)$. Then $a$ lifts to an invariant absolute connection
$b\in \H^1(G, \sO^*_G \to \OG)^{{\rm inv}}$. 
\end{lem}
\begin{proof}  One has as in \eqref{2.19} 
\ml{(\mu^*-p_1^*-p_2^*)(b')\in \ 
{\rm Im \ }H^0(\sO_{G\times G})\otimes \OK  \\
\subset
\H^1(G\times G, \sO^*_{G\times G} \to \Omega^1_{G\times G}) 
}{2.37}
Now 
\begin{gather}
{\rm Im \ }H^0(\sO_{G\times G})\otimes \OK = \notag \\
H^0(\sO_{G\times G})\otimes \OK/d\log {\rm Ker \ }
\{H^0(\sO_{G\times G}^*) \to H^0(\Omega^1_{G \times G/K})\}= \notag \\
H^0(\sO_{G\times G})\otimes (\OK/d\log K^*). \notag
\end{gather}
Exactness of the sequence in \eqref{2.19} implies that
there exists an element $x \in H^0(G, \sO\otimes \OK)$ with
$(\mu^*-p_1^*-p_2^*)(x)=(\mu^*-p_1^*-p_2^*)(b')$. Take $b=b'-x$. 
\end{proof}
\begin{prop}\label{prop3} One has an exact sequence
\begin{multline}\label{2.38} 0 \to H^0(\sO_G)^{{\rm inv}}\otimes \OK \to
\H^1(G,\sO_G^* \to \OG)^{{\rm inv}} \\ 
\to \H^1(G,\sO_G^* \to \OGK)^{{\rm inv}} \to H^1(\sO_G)^{{\rm inv}}\otimes \OK
\end{multline}
\end{prop}
\begin{proof} This is immediate from the lemma. 
\end{proof}

Recall our notation. $C$ is a smooth, projective, geometrically
connected curve over $K$.  $G=J_\sD$ with $\sD=\sum m_i c_i$, $D=\sum c_i$,
$\sD' = \sD - D$.  

\begin{prop}\label{prop2.12} There exists a diagram of exact sequences,
with vertical 
arrows isomorphisms: 
\begin{tiny}
\eq{\begin{array}{cccccc}
0 & \to & H^0(\sO_G)^{{\rm inv}} \otimes (\OK/d \log K^*) 
& \to & \H^1(G,\sO_G^* \to
\OG)^{{\rm inv}} \\
&& \rule{0cm}{.5cm}\downarrow e  && \downarrow  \\
0 & \to & \Big(H^0(\sO_C(\sD'))/H^0(\sO_C)\Big) \otimes (\OK/d \log K^*) & \to &
\rule{0cm}{.5cm}\makebox[3.5cm][c]{$\frac{\H^1(C,j_*(\sO_{C-D}^*) \to
\Omega^1_C\ldr(\sD'))}{(\OK/K^\times)}$} \\
\rule{0cm}{1cm} \\
& \to & \H^1(G,\sO_G^* \to \OGK)^{{\rm inv}} & \to &
H^1(\sO_G)^{{\rm inv}}\otimes \OK  \\
&& \rule{0cm}{.5cm}\downarrow  && \downarrow  h\\
&\to & \rule{0cm}{.5cm}\H^1\Big(C,j_*(\sO_{C-D}^*)\to
\Omega^1_{C/K}(\sD)\Big) & \to & 
H^1(\sO_C(\sD'))\otimes \OK \\
\end{array}
}{2.39}
\end{tiny}
\end{prop}

\begin{proof} The first step is to compute $H^i(\sO_G)^{{\rm inv}}$ for
$i=0,1$. Let $W$ be a finite dimensional
$K$-vector space, and suppose $G=J_\sD$ is a vectorial extension 
\eq{0 \to W\otimes \G_a \to G \stackrel{p}{\to} G_0 \to 0
}{2.40}
We know by lemma \ref{lem2.1} that this sequence pulls back from
an extension of $J$ by $W\otimes \G_a$. Let
\eq{0 \to \sO_J \to fil_1 \to W^*\otimes \sO_J \to 0
}{2.41}
be the exact sequence of functions of filtration degree $\le 1$ as in
lemma \ref{lem2.3}, and let $\p : W^* \to H^1(\sO_J)$ be the
boundary map in cohomology. 
\begin{lem}We have
\begin{gather}H^0(\sO_G)^{{\rm inv}} \cong 
\label{2.42}\ker(\partial: W^* \to H^1(\sO_J)) \\
\label{2.43} H^1(\sO_G)^{{\rm inv}} \cong H^1(\sO_J)/\partial(W^*). 
\end{gather}
\end{lem}
\begin{proof}[proof of lemma]
One has a filtration $fil_\cdot p_*\sO_G$ with $fil_0 = \sO_{G_0}$ and
$gr_r = \sym^r(W^*)\otimes \sO_{G_0}$. The corresponding spectral
sequence looks like 
\begin{equation}\label{2.44} E_1^{pq} = H^{p+q}(G_0, gr_{-p}) =
H^{p+q}(G_0, \sym^{-p}(W^*)\otimes \sO_{G_0}) \Rightarrow
H^{p+q}(\sO_G). 
\end{equation}
Let 
\eq{0 \to H_0 \to G_0 \to S \to 0
}{2.45}
be as in lemma \ref{lem2.2}, so $S$ is the maximal quotient torus
of $G_0$. 

The equation  \eqref{2.6} identifies
$H^i(G_0, \sO_{G_0})$ with $H^i(J \times S, \sO_{J \times S})$,
and the invariance condition might be looked at on $J \times S$. 
Let us write $H^0(S, \sO_S) = K \oplus V, f \mapsto f(0) \oplus (f - f(0))$,
where $V$ consists of the regular functions which vanish 
at $\{0\} \in S$. 
Then $H^i(G_0, \sO_{G_0})^{{\rm inv}}= H^i(J, \sO_J)^{{\rm inv}} 
\oplus (H^i(J, \sO_J) \otimes V)^{{\rm inv}}$.  
Thus if a class $F= \sum \varphi_f \otimes f \in H^i(\sO_J)
\otimes V$ is invariant, where $\varphi_f \in H^i(J, \sO_J)$
and the $f \in V$ are linearly independent over $K$, 
then 
\begin{gather}
(\mu^* - p_1^* -p_2^*)(F)|(J\times J \times \{0\} \times S)
= \notag \\
\sum (\mu^* - p_2^*)(\varphi_f) \otimes f = 0 \notag 
\end{gather}
thus $\mu^*\varphi_f=p_i^*\varphi_f$ and
$\mu^*\varphi_f|\{0\}\times J= \varphi_f = \varphi_f|\{0\}.$
So for $i\geq 1$, this implies that $\varphi_f=0$ and for $i=0$
this implies that $F\in H^0(\sO_S)^{{\rm inv}}=0.$ 
  
In short: 
\eq{H^*(\sO_{G_0})^{{\rm inv}} \cong \Big(H^*(\sO_J)\otimes
H^0(\sO_S)\Big)^{{\rm inv}} \cong H^*(\sO_{H_0})^{{\rm inv}}
\cong H^*(\sO_J)^{{\rm inv}}
}{2.46}
Thus, it suffices to prove the lemma with $G_0$ replaced by
$H_0$, so we may assume the quotient torus $S=(0)$.
Since in this case $H^i(\sO_{G_0})\cong \wedge^i H^1(\sO_J)$, one sees
that 
 pullback under the multiplication by $N$ map, $N\delta:G \to G$ acts
 on $E_1^{pq}$ by multiplication by $N^q$. It follows that the
 spectral sequence \eqref{2.44} degenerates at $E_2$. In particular the
 eigenspace where $N\delta^*$ 
 acts by multiplication by $N$ on $H^1(\sO_G)$ is
\begin{equation}\label{2.47} H^1(\sO_J)/\partial(W^*) \cong E_2^{0,1}
 \cong E_\infty^{0,1} \hookrightarrow H^1(\sO_G). 
\end{equation}

Note that as a quotient of $H^1(\sO_J)$ the space $E_\infty^{0,1}$ is
clearly invariant. Conversely, let $\Delta: G \to G\times G$ be the
diagonal.  Since $\mu\circ\Delta = 2\delta$ it follows that for $a\in
H^1(\sO_G)^{{\rm inv}}$ we have 
\eq{(2\delta)^*(a) = \Delta^*\mu^*(a) = \Delta^*(p_1^*(a)+p_2^*(a)) = 2a
}{2.48}
so necessarily $a\in E_\infty^{0,1}$, proving \eqref{2.43}. A similar
argument on $E_\infty^{-1,1}$ proves \eqref{2.42}.

We remark here that $H^0(J, fil_1)$ is in a natural way a subspace of
the regular functions on $G$, and \eqref{2.42} takes the quotient
of this by $H^0(J, \sO_J)= K$. This is because we have forced the rigidification 
condition in the definition \ref{def2.6}.
\end{proof}

We return to the proof of proposition \ref{prop2.12}. The exact
sequence
\begin{equation}\label{2.49}
0 \to \sO_C \to \sO_C(\sD') \to \sO_C(\sD')/\sO_C \to 0
\end{equation}
defines a map
\begin{equation}\label{2.50}
\psi: W^* := H^0(\sO_C(\sD')/\sO_C ) \to H^1(\sO_C).
\end{equation}
By lemma \ref{lem2.1}, as a group extension of $G_0$,
the group $G$ is defined by a unique map 
from $H^0(C, \Omega^1_{C/K})$ to a vector space. We claim that
this map is the dual of $\psi$. To see this, one identifies $J$
and $J^\vee$. Then it is well know that the universal
vectorextension on $J^\vee$ is
$$0 \to H^0(C, \Omega^1_{C/K}) \to \H^1(C, \sO^*_C \to
\Omega^1_{C/K}) \to {\rm Pic}^0(C) \to 0$$
inducing the universal vectorextension 
$$0 \to H^0(C, \Omega^1_{C/K}) \to \H^1(C, \sO^*_{C,D} \to
\Omega^1_{C/K}) \to {\rm Pic}^0(C, D) \to 0$$
on 
$${\rm Pic}^0(C, D):= {\rm Ker}\Big(H^1(C, \sO^*_{C,D}) 
\to H^1(C, \Omega^1_{C/K})\Big)$$ 
where $$\sO^*_{C, Z}:= {\rm Ker} (\sO^*_C \to \sO^*_{Z})$$
for any subscheme $Z \subset C$.
The  map of complexes
$$a: \{\sO^*_{C, \sD} \to 0\} \to \{ \sO^*_{C, D} \to
\Omega^1_{C/K}|_{\sD'}\}$$
induces an isomorphism on $\H^1$. Indeed, $a$ sends
the exact sequence
$$0 \to H^0(C, (1+\sO_{\sD'}(-D))) \to H^1(C, \sO^*_{C, \sD})
\to H^1(C, \sO^*_{C,D}) \to 0$$
to the exact sequence
$$0 \to H^0(C, \Omega^1_{C/K}|_{\sD'}) \to 
\H^1(C,  \sO^*_{C, D} \to
\Omega^1_{C/K}|_{\sD'}) \to H^1(C, \sO^*_{C,D}) \to 0,$$
so one has just to see that
$$ d \log: H^0(C, (1+\sO_{\sD'}(-D))) \to
H^0(C, \Omega^1_{C/K}|_{\sD'})$$
is an isomorphism.
But $ H^0(C, \sO_{\sD'}(-D)) \cong H^0(C, (1+\sO_{\sD'}(-D)))$ via
the exponential map and the quasiisomorphism
\cite{DeI}
$$\{\sO_C(-\sD)\to \Omega^1_{C/K}(-\sD')\} \to \{
\sO_C(-D) \to \Omega^1_{C/K}\}$$
allows to conclude.

Define a bundle $\sE$ on $C$ by pullback
\begin{equation}\label{2.51}\begin{CD} 
0 @>>> \sO_C @>>> \sE @>>> W^*\otimes\sO_C @>>> 0 \\
@. @| @VVV @VVV \\
0 @>>> \sO_C @>>> \sO(\sD') @>>> \sO_C(\sD')/\sO_C @>>> 0
\end{CD}
\end{equation}
Because of the isomorphism $H^1(\sO_J)\cong H^1(\sO_C)$, the top row
of the diagram \eqref{2.51} pulls back uniquely from an extension of
$W^*\otimes \sO_J$ by $\sO_J$. There is a unique vectorial extension 
\begin{equation}\label{2.52}
0 \to W\otimes \G_a \to H \stackrel{r}{\to} J \to 0
\end{equation}
such that the above extension of vector bundles coincides with 
\begin{equation}\label{2.53}
0 \to \sO_J \to fil_1r_*\sO_H \to W^*\otimes \sO_J \to 0
\end{equation}
{From} this we get a diagram (defining $t$ and $u$. Here $i:C
\hookrightarrow J$) 
\begin{equation}\label{2.54}
\begin{CD} 0 @>>> \sO_J @>>> fil_1r_*\sO_H @>>> W^* \otimes \sO_J @>>>
0 \\ 
@. @VVV @VV t V @VV u V \\
0 @>>> \sO_C @>>> \sO_C(\sD') @>>> \sO_C(\sD')/\sO_C @>>> 0
\end{CD}
\end{equation}

We get a diagram with exact rows 
\begin{tiny}
\begin{equation}\label{2.55}
\begin{array}{ccccccccccccc} 0 & \to & H^0(\sO_{J}) & \to &
H^0(fil_1r_*\sO_H) & \to & H^0(W^*\otimes \sO_{J}) &
\stackrel{\partial}{\to}  & H^1(\sO_{J}) & \to  \\
&&\downarrow \cong && \downarrow t && \downarrow \cong && \downarrow
\cong  \\ 
0 & \to & H^0(\sO_C) & \to & H^0(\sO_C(\sD')) & \to &
H^0(\sO_C(\sD')/\sO_C) & \to & H^1(\sO_C) & \to \\
\\
 H^1(\sO_H)^{{\rm inv}} & \to & 0 \\
\downarrow v \\ 
 H^1(\sO_C(\sD')) & \to & 0  
\end{array}
\end{equation}
\end{tiny}
The diagram \eqref{2.55} gives isomorphisms
\begin{gather}\label{2.56}e: H^0(\sO_G)^{{\rm inv}} \cong \ker(\partial)
\cong {\rm Ker}(H^0(\sO_C(\sD')/\sO_C)\to H^1(\sO_C)) \\  
\label{2.57} h: H^1(\sO_G)^{{\rm inv}} \cong H^1(\sO_C(\sD')). 
\end{gather}
These are two of the desired arrows for the diagram in the
proposition.

\begin{lem}The natural map on relative connections
\eq{\H^1(G, \sO_G^* \to \OGK)^{{\rm inv}} \to \H^1\Big(C,j_*(\sO_{C-D}^*) \to
\Omega^1_{C/K}(\sD)\Big) 
}{2.58} 
is an isomorphism. 
\end{lem}
\begin{proof}[proof of lemma] We note the following facts:
\begin{enumerate}\item $H^1(G, \OGK)^{{\rm inv}}=\Big( H^0(G, \OGK)^{{\rm inv}}\otimes
H^1(\sO_G) \Big)^{{\rm inv}}= (0)$.
\item $H^1(\sO_G^*)^{{\rm inv}}\cong H^1(\sO_{C-D}^*)$. Indeed, as remarked
in the proof of lemma \ref{lem2.8} one has $J(K) \twoheadrightarrow
H^1(\sO_G^*)^{{\rm inv}}$. One checks that the kernel is generated by
divisors of degree $0$ supported on $D$.  
\item $H^0(\sO_G^*)/\text{consts.}\cong H^0(\sO_G^*)^{{\rm inv}}\cong
H^0(\sO_{C-D}^*)/\text{consts.}$. 
\item 
$$\Big(H^0(\OGK)/d\log(H^0(\sO_G^*))\Big)^{{\rm
inv}}=H^0(\OGK)^{{\rm inv}}/d\log 
(H^0(\sO_G^*)^{{\rm inv}}) 
$$
(This is seen by noting $H^0(\sO_G^*)^{{\rm inv}}= 
{\rm Hom}(G,\G_m)$, so one has
a homomorphism $\psi : G \to \G_m^r$ such that $\psi^*$ is an
isomorphism on global units modulo constants. The assertion then
reduces to the case 
$G=\G_m^r$, which is easy.)  
\end{enumerate} 

We build a diagram \begin{tiny}
\begin{equation}\label{2.59} \begin{array}{ccccccccc}
 0 & \to & \frac{H^0(\OGK)^{{\rm inv}}}{H^0(\sO_G^*)^{{\rm inv}}} & \to &
\H^1(\sO_G^* \to \OGK)^{{\rm inv}} & \to & H^1(\sO_G^*)^{{\rm
inv}} & \to & 0 \\ 
&& \downarrow && \downarrow && \downarrow \\
0 & \to & \frac{H^0(\Omega^1_{C/K}(\sD))}{H^0(\sO_{C-D}^*)} & \to &
\H^1\Big(C,j_*(\sO_{C-D}^*) \to \Omega^1_{C/K}(\sD)\Big) & \to &
H^1(C-D,\sO^*) & \to &  0 
\end{array}\end{equation}
\end{tiny}
Since the left and right hand arrows are isomorphisms, it follows that
the central arrow is as well, proving the lemma. 
\end{proof}
The assertions of the proposition follow easily from the lemma. 
\end{proof}

Finally, we need an analogous result for integrable connections. More
precisely, we consider a slightly weaker condition. 
\begin{defn}\label{def2.15} Let $X$ be a variety over $K$. A
connection $\nabla : E \to E\otimes \Omega^1_X\otimes K(X)$ (so
possibly with poles) is said to have
vertical curvature if the curvature  
\eq{\nabla^2 : E \to 
E\otimes\Omega^2_X\otimes K(X)
}{2.60}
has values in the subsheaf $E\otimes\Omega^2_K\otimes K(X) \subset 
E\otimes\Omega^2_X\otimes K(X)$.
The group of line bundles with vertical curvature will be denoted
$$\H^1(X,\sO_X^* \to \Omega^1_X)^{{\rm vert \ }}$$
and similarly for invariant line bundles with vertical
curvature
$$\H^1(X,\sO_G^* \to \Omega^1_G)^{{\rm inv,vert \ }}.$$ 
\end{defn}
\begin{prop}\label{prop2.16} 
With notation as above, we have isomorphisms
\begin{gather}\label{2.61} \H^1(G, \sO^*_G
\to \OGK)^{{\rm inv,vert}}= \H^1(C, j_*\sO^*_{C-D}
\to \Omega^1_{C/K}(\sD))^{{\rm vert}} \\
\label{2.62} \H^1(G, \sO^*_G \to \OG)^{{\rm inv}, {\rm vert \ }}\cong
\frac{\H^1(C,j_*\sO_{C-D}^* \to 
\Omega^1_C\ldr(\sD'))^{{\rm vert \ }}}{\OK/K^\times}
\end{gather}
\end{prop}
\begin{proof}For example, in the absolute case, the curvature of a
line bundle with invariant absolute connection on $G$ is a 
section $\eta \in H^0(\Omega^2_G/\Omega^2_K\otimes\sO_G)$ satisfying
$\mu^*(\eta) = p^*_1(\eta)+p^*_2(\eta)$. It is easy to see that such a
section lies in the subsheaf $\OGK\otimes\OK$. The isomorphism
\eqref{2.62} follows from proposition \ref{prop2.12} and the fact that
pullback to $C$ of invariant forms is injective by \eqref{2.32}. The
case of relative forms is similar and is left for the reader. 
\end{proof}

\section{The Geometric Setup}

We continue to work with a curve $C/K$ and a line bundle $L$ on $C$
of degree $0$. Let $\nabla_{/K}:L \to L\otimes\OCK(\sD)$, where $\sD =
\sum m_ic_i$. As in the previous section, write $D=\sum c_i$ and $\sD'
= \sD - D$. 
\begin{lem}\label{lem3.1} Assume $\nabla_{/K}|C-D$ lifts to an
absolute integrable connection $\nabla' : L_{C-D} \to
L_{C-D}\otimes\Omega^1_{C-D/k}$. Then $\nabla'$ extends to an
absolute integrable connection
\eq{\nabla : L \to L\otimes\Omega^1_{C/k}\ldr(\sD').
}{3.1}
The notation $\ldr$ refers to log poles at $D$ as in \eqref{1.4}. 
\end{lem}
\begin{proof}Let $e$ be a basis for $L$ at $c$ a point with
multiplicity $m\ge 1$ in $\sD$, and let $x$ be a
local parameter at $c$ on $C$. Write
\eq{\nabla'(e) = A(x)dx+\sum_iB_i(x)d\tau_i;\quad d\tau_i \text{ basis
in } \OK,\
x^mA(x)\in \sO_{C,c}. 
}{3.2}
We must show $x^{m-1}B_i(x)$ is regular at $c$. But integrability of
$\nabla'$ implies that $\p A/\p \tau_i = \p B_i/\p x$, from which the
assertion is clear. 
\end{proof} 

We know from proposition \ref{prop2.16} 
that the restriction to $C-D$ of an integrable
absolute connection of the form 
\eqref{3.1} pulls back from a unique invariant integrable absolute
connection $\sL \to \sL\otimes\OG$ on $G=J_\sD$. More precisely, we fix a
basepoint $c_0 \in (C-D)(K)$ and normalize our connection \eqref{3.1}
to be trivial at the basepoint by tensoring with a pullback from
$\Spec(K)$. 

We consider now the basic geometric picture of Deligne \cite{De}
\eq{\pi:\sym^N(C-D) \to G_N;\quad N=2g-2+\sum_im_i
}{3.3}
where $G_N$ is the $J_\sD$-torseur of degree $N$ line bundles
trivialized along $\sD$ and $\pi(\sum z_i)=\sO(\sum z_i)$ with
trivialization given by restricting to $\sD$ the canonical (upto
scalar in $K$) map $\sO_C \to \sO_C(\sum z_i)$. Note that $N=\deg(\OCK(\sD))$
and $\dim G = g-1+\sum m_i = N-g+1$. (Recall we assume $\sD\neq
\emptyset$.) We identify $G_N \cong G$ by sending the point
$[\sO(Nc_0)]\mapsto 0$, and we write $\sL$ for the resulting line
bundle with connection on $G_N$. The basic remark of Deligne is
\begin{prop}\label{prop3.2}Assume
$$H^0_{DR}(C-D,(L,\nabla))=H^2_{DR}(C-D,(L,\nabla))  =
0.$$ 
Then 
\ml{\det(H^*_{DR}(C-D,(L,\nabla))=\det(H^1_{DR}((C-D,(L,\nabla)))) \\
\cong H^N(\sym^N(C-D),(\pi^*(\sL, \nabla)))
}{3.4}
as a line with connection on $K$.
\end{prop}
\begin{proof} Our hypotheses imply $H^1_{DR}(C-D,(L,\nabla))$ has
dimension $N$. Consider the diagram
\eq{\begin{array}{ccc} (C-D)^N & \stackrel{p}{\to} &
\makebox[2.4cm][r]{$\sym^N(C-D)$} \\ 
&  \searrow q & \makebox[2.4cm][l]{$\downarrow \pi$} \\
&&  \makebox[2.4cm][l]{$G_N$} 
\end{array}
}{3.5}
We have $q^*(\sL) \cong L\boxtimes\cdots\boxtimes L$ (exterior tensor
product on $(C-D)^N$. The K\"unneth formula gives
\eq{H^N((C-D)^N,(L\boxtimes\cdots\boxtimes L,\nabla)) \cong
H^1_{DR}(C-D,(L,\nabla))^{\otimes N}. 
}{3.6}
There is an action of the symmetric group $\sS_N$ on the pair
$$((C-D)^N, L\boxtimes\cdots\boxtimes L).$$ The resulting action on
$(H^1_{DR})^{\otimes N}$ is alternating because of the odd degree
cohomology, so the invariants are precisely $\det H^1_{DR}$. There is
an evident map 
\ml{p^*: H^N_{DR}(\sym^N(C-D),\pi^*\sL) \to
H^N_{DR}((C-D)^N,L\boxtimes\cdots\boxtimes L)^{\sS_N} \\
= \det H^1_{DR}(C-D,L)
}{3.7}
To show this map is an isomorphism, it suffices to remark that one has
a trace map 
\eq{p_* :  H^N_{DR}((C-D)^N,L\boxtimes\cdots\boxtimes L) \to
H^N_{DR}(\sym^N(C-D),\pi^*\sL)  
}{3.8}
Because $L\boxtimes\cdots\boxtimes L = p^*\pi^*\sL$, the existence of
such a trace follows from the projection formula and the trace in de
Rham cohomology with constant coefficients. 
\end{proof}

Now one uses the geometry of the map $\pi$ and \eqref{3.4} to compute
the determinant.

\begin{lem}\label{lem3.3} Let $X$ be a smooth variety over a field of
characteristic $0$. Let $A\subset X$ be a smooth subvariety of
codimension $p$. Let $(E,\nabla)$ be an integrable connection on
$X$. Then 
\eq{\H^n_A(X, E\otimes \Omega^*_X) \cong H^{n-2p}(A, E\otimes
\Omega^*_A). 
}{3.9}
\end{lem}
\begin{proof} Write $\uh_A^r(F)$ for the Zariski sheaf associated to the
presheaf $U \mapsto H^r_A(U,F)$ for any Zariski sheaf $F$ on $X$. For
$F$ locally free, $\uh_A^r(F)= (0)$ for $r\neq p$ by purity. Duality
theory gives  (here $X\supset A_\alpha \supset A$
runs through 
nilpotent thickenings) 
\ml{E\otimes \Omega^m_A \to {\rm Ext}^p(\sO_A, E\otimes \Omega^{m+p}_X)
\to \varinjlim_\alpha {\rm Ext}^p(\sO_{A_\alpha},E\otimes\Omega^{m+p}_X) \\
\cong \uh_A^p(E\otimes\Omega^{m+p}_X).
}{3.10}
We want to show that this map is an isomorphism, compatible
with the connection, thus yielding a quasiisomorphism
of complexes
\eq{E\otimes\Omega^*_A \to \uh_A^p(E\otimes\Omega^*_X).
}{3.11}
The problem is local, so we can assume
$$A\subset A_1 \subset \ldots \subset A_p = X
$$ 
with $A_i$ smooth of codimension $p-i$ in $X$. Now
$\uh_A^p(E\otimes\Omega^*_X)$ represents
$R\underline{\underline{\Gamma}}_A(E\otimes\Omega^*_X)[p]$ in the
derived category of Zariski sheaves on $A$, and in the derived
category we may write
\eq{R\underline{\underline{\Gamma}}_A(E\otimes\Omega^*_X)[p] =
R\underline{\underline{\Gamma}}_A[1]\circ
R\underline{\underline{\Gamma}}_{A_1} [1]\circ\ldots\circ
R\underline{\underline{\Gamma}}_{A_{p-1}}(E\otimes\Omega^*_X)[1]
}{3.12}
In this way, we reduce to verifying \eqref{3.9} in the case $p=1$. So,
suppose $A:t=0$ in $X=\Spec(R)$. We have
\eq{ \uh_A^1(E\otimes\Omega^*_X)\cong
E_{R[t^{-1}]}\otimes
\Omega^*_{R[t^{-1}]}/E_R\otimes\Omega^*_R 
}{3.13}  
as $ H^1_{DR}(X, E) \subset H^1_{DR}(X-A, E).$

By \cite{DeI},
since $E$ has no singularity along $t=0$, one has
\eq{(E_{R}\otimes \Omega^*_{R}(\log(t=0)), \nabla) 
\stackrel{{\rm q.iso.}}{\to}
(E_{R[t^{-1}]}\otimes \Omega^*_{R[t^{-1}]}, \nabla|_{{\rm Spec}R[t^{-1}]}) 
}{3.14} 
Thus 
${\rm res \ }: \uh_A^1(E\otimes\Omega^*_X, \nabla)) \to 
E_{R/tR} \otimes (\Omega^*_{R/tR}[-1], \nabla|_{{\rm Spec}(R/tR)}) 
$
is an isomorphism.
\end{proof}

\begin{lem}\label{lem3.4}
Let $p:G_N \to J_N$ be the projection to the corresponding
torseur over the absolute jacobian. Write $[\OCK(\sD)]\in J_N$ for the
point corresponding to the canonical bundle twisted by $\sO(\sD)$. Let
$a\in G_N$. We have
\eq{\pi^{-1}(a) = \begin{cases} \A^{g-1} & p(a) \neq [\OCK(\sD)] \\
\A^g & p(a) = [\OCK(\sD)];\ \p(a) = 0 \\
\emptyset &  p(a) = [\OCK(\sD)];\ \p(a) \neq 0
\end{cases}
}{3.15}
where by definition $\A^{g-1}= \emptyset$ is $g=0$.
Note that if $p(a)=[\OCK(\sD)]$, then $a$ corresponds to a
trivialization $\sO_\sD \cong \OCK(\sD)|_\sD$ defined upto
scalar. $\p(a)$ in the above refers to the evident boundary of this
trivialization in $H^1(\OCK)=K$ (again upto scale).
\end{lem}
\begin{proof}Let $p(a)$ correspond to a line bundle $M$ of degree $N$,
we consider the exact sequence
\eq{0 \to M(-\sD) \to M \to M|_\sD \to 0
}{3.16}
Suppose first $M\neq \OCK(\sD)$. Then $H^1(M(-\sD))=(0)$, so any
trivialization in $H^0(M|_\sD)$ lifts to $H^0(M)$, and the space of
such liftings is a torseur under $H^0(M(-\sD))$, a vector space of
dimension $g-1$. (Note this is an affine torseur, not a projective
torseur.) If $M = \OCK(\sD)$, $H^0(M(-\sD))$ has dimension $g$, and
the image $H^0(M) \to H^0(M|_\sD)$ has codimension $1$. 
\end{proof}

\begin{remark}\label{rmk3.5} If we choose local parameters $t_i$ at
$c_i \in \sD$, then 
$$H^0(\OCK(\sD)|_\sD)$$ 
can be identified with the
space of polar parts of $1$-forms with poles along $\sD$, and the map
$\p$ is given by the residue
\eq{\p(\sum_i\sum_{j=0}^{m_i-1} u_{ij}dt_i/t_i^{m_i-j}) = \sum_i
u_{i,m_i-1} 
}{3.17}
Note the (open) condition for an element in $H^0(\OCK(\sD)|_\sD)$ to
be a trivialization is simply
\eq{\prod_i u_{i0} \neq 0.
}{3.18}
Because $B\subset G_N$, we must factor out by the action of $\G_m$,
which we can normalize away by setting $u_{10} = 1$. 
Thus we have
\ml{B:= \pi(\sym^N(C-D))\cap p^{-1}[\OCK(\sD)] = \\
\Big\{\sum_i\sum_{j=0}^{m_i-1} u_{ij}dt_i/t_i^{m_i-j} \ \Big| \
\sum_i u_{im_i-1}= 0;\ 
\prod_i u_{i0} \neq 0;\ u_{10}=1 \Big\}.
}{3.19}
\end{remark}

Define $A:= \pi^{-1}(B)\subset \sym^N(C-D)=:X$. Using the localization
sequence and lemma \ref{lem3.4}, we get an exact sequence
\eq{\minCDarrowwidth.5cm \begin{CD}
H^{N-2g}_{DR}(A/K,\pi^*\sL) @>>> H^N_{DR}(X/K,\pi^*\sL) @>>>
H^N_{DR}(X - A/K,\pi^*\sL) \\
@AA\cong A @. @AA\cong A \\
H^{N-2g}_{DR}(B/K,\sL|_B) @.@.  \makebox[3cm][r]{$H^N_{DR}(G_N -
p^{-1}[\OCK(\sD)]/K,\sL)$} 
\end{CD}
}{3.20} 
where $p:G_N \to J_N$ with $J= J(C)$ the absolute Jacobian. The fact
that the vertical arrows in this diagram are isomorphisms follows
because the maps are maps of affine space bundles and the line bundles
with connection are pulled back from the base.  

To simplify the presentation, we will assume that $m_1\ge 2$. Another
proof of our formula for the de Rham determinant in the case $\sD = D
= \sum c_i$, (i.e. for regular singular points) will be given in
theorem \ref{thm4.5}. 

\begin{lem}\label{lem3.6} Assume the line bundle $L$ on $C$ has degree
$0$, and that $\sD$ is minimal, i.e. $\nabla: L|_\sD \cong
L\otimes\OCK(\sD)|_\sD$. We continue to assume also that $\sD=\sum
m_ic_i$ with $m_1\ge 2$. Then
\eq{H^*_{DR}(\sym^N(C-D)-A/K,\pi^*\sL) = (0).
}{3.21}
\end{lem}
\begin{proof}The isomorphism on the right in \eqref{3.20} implies we
must show $H^*_{DR}(G_N - p^{-1}[\OCK(\sD)]/K,\sL) = (0)$. The
assumption $m_1\ge 2$ means we have a $\G_a$ action by translation on
$ G_N - p^{-1}[\OCK(\sD)]/K,\sL $, and minimality of $\sD$ implies
that the connection is nontrivial on the fibres. The fibration is 
Zariski-locally trivial, so the Leray spectral sequence for de Rham
cohomology reduces us to showing $H^*_{DR}(\G_{a,S}/S,(\sO,\Xi))=(0)$
where $\Xi$ is an everywhere non-zero, translation invariant, relative
$1$-form on 
$\sO_{\G_{a,S}}$.  In other words, for $s\in \sO_S^\times$ we must show
\eq{\sO_S[t] \stackrel{d+sdt}{\longrightarrow} \sO_S[t]dt
}{3.22}
has trivial cohomology. This is straightforward.
\end{proof}

\begin{lem}\label{lem3.7} Assume $H^0_{DR}(C-D,L)=(0)$. Then  
\eq{H^{m-2}_{DR}(B/K,\sL|_B) \cong
H^N_{DR}(\sym^N(C-D)/K,\pi^*\sL)\cong K 
}{3.23}
as a line with a connection over $K$.
\end{lem}
\begin{proof}Note $m-2=N-2g$. Extending the top sequence in
\eqref{3.20} one step to 
the left and using the previous lemma gives the left isomorphism. We
have already seen the isomorphism on the right. 
\end{proof}

Our task now is to calculate $H^{m-2}_{DR}(B/K,\sL|_B)$ with its
connection. We assume that $\sL \in {\rm Pic}^0(J)$ as in lemma
\ref{lem3.6}. Then $\sL$ carries a relative invariant connection 
$d_{/K}$ on $J$, and 
$\nabla_{/K}= d_{/K}+ \Xi$ for some invariant form 
$\Xi \in H^0(G, \Omega^1_{G/K})^{{\rm inv}}$.
Changing the choice of $d_{/K}$ changes $\Xi$
to $\Xi + p^*(\alpha)$, where $\alpha \in H^0(J, \Omega^1_{J/K})^{{\rm
inv}}$ and $p: G_N \to J_N$ is the torseur under the affine
group $\sG:= \ker(G \to J)$. 
In particular, $\Xi|p^{-1}[\Omega^1_{C/K}(\sD)]$ does not
depend on the choice of $d_{/K}$. 
As $p^{-1}(\OGK(\sD))$ is isomorphic to $\sG$, and 
we see that
\eq{(\sL,\nabla_{/K})|_B = (\sO_B, d+ \Xi|_B).
}{3.24} 
We say that $\Xi|_B$ vanishes at a point $b\in B$ if $\Xi(b)=0$ in
$\mathfrak
m_b/\mathfrak m_b^2$. 
\begin{lem}\label{lem3.8} Let $b\in B\subset G_N$ correspond to the
trivialization on $\OCK(\sD)|_\sD$ given by $\nabla|_\sD : L|_\sD \to
L\otimes \OCK(\sD)|_\sD$. Then $\Xi|_B$ vanishes at $b$. $\Xi|_B$ does
not vanish at any other point of $B$. 
\end{lem}
\begin{proof}
Let us write $\eta:=i^*\Xi|C-D$, where $i: C-D \to G$ is the cycle 
map. Then by definition, the trivialization
of $\OCK(\sD)|_\sD$
associated to
$i^*(\sL,\nabla)$ depends only on $\eta|\sD$, or equivalently
only on $\Xi|p^{-1}[\Omega^1_{C/K}(\sD)]$. 
We have
\eq{\pi^*\Xi = \sum_{i=1}^N \eta_i
}{3.25}
where $\eta_i$ is the pullback of $\eta$ via the $i$-th projection
$(C-D)^N \to C-D$. 

Suppose for a moment that the divisor of $\eta$ (viewed as a section
of $\OCK(\sD)$) is reduced, $(\eta) = \sum e_i;\ e_i \in
(C-D)(\bar{K})$. Let $e := (e_1,\cdots,e_N)\in \sym^N(C-D)$ be the
point corresponding to $\eta$. We have $e\in A=\pi^{-1}(B) \subset
\sym^N(C-D) $ and $b=\pi(e)$. Since $A \to B$ is a projective
bundle, there is a 
surjection on tangent spaces 
\eq{\begin{array}{ccc}T_A(e) & \inj & T_{\text{Sym}^N(C-D)}(e) \\
\makebox[0cm][r]{surjective}\downarrow && \downarrow \pi_*  \\ 
T_B(b)  & \inj & T_{G_N}(b) \end{array}
}{3.26}
Since $T_{\sym^N(C-D)}(e)$ is spanned by expressions $\sum
\tau_i|_{e_i}$, to show $\Xi|_B$ vanishes at $b$, it suffices to show
\eq{<\Xi,\pi_*(\sum \tau_i|_{e_i})> = 0.
}{3.27} 
This expression equals
\eq{<\pi^*(\Xi),\sum \tau_i|_{e_i}> = \sum <\eta,\tau_i|_{e_i}>.
}{3.28}
Each term on the right vanishes because $\eta(e_i)=0$ in $\mathfrak
m_{e_i}/\mathfrak m_{e_i}^2$. The general case ($(\eta)$ not necessarily
reduced) follows from this by a specialization argument.

We postpone until lemma \ref{lem3.10} 
the proof that $\Xi|_B$ doesn't vanish at any
other point of $B$. 
\end{proof}

By assumption we start with an absolute, invariant, integrable
connection on $\sL$ of degree 0 on $J$. Restricting to $B$, we get an
absolute closed
invariant $1$-form $\Psi$, whose corresponding relative form is
$\Xi$. Recall \eqref{3.19} we have coordinates $u_{ij}$ on $B$ with
$u_{10}=1$, $\prod u_{i0}\neq 0$, and $u_{1,m_1-1}=-\sum_{i\ge
2}u_{i,m_i-1}$. 
\begin{lem}\label{lem3.9}
Under the assumption that 
$$\{i^*\nabla_{/K}: i^*\sL \to i^*\sL \otimes
\Omega^1_{C/K}(\sD)\} \to \{j_*i^*\sL|(C-D) 
\to j_*i^*\sL \otimes \Omega^1_{(C-D)/K}\}$$
is a quasiisomorphism,  and $i^*\nabla_{/K}$ has poles
along all points of $D$,
we can arrange that a $K$-basis for 
$$H^{m-2}_{DR}(B/K,\sL|_B) = H^{m-2}_{DR}(B/K,(\sO_B,d+\Xi))
$$
is given by the closed form
\eq{\theta:= \prod_{\substack{(i,j) \\ m_i\ge 2}}
du_{ij}\wedge\prod_{\substack{i \\ m_i=1}} 
\frac{du_{i0}}{u_{i0}} 
}{3.29}
\end{lem}
\begin{proof} Recall \eqref{2.29} the relative invariant forms on
$\sG:=\ker(G \to J)$ are the $\tau_{ij}$ defined by
the expression
\eq{\sum_{j=0}^{m_i-1} \tau_{ij}T_i^j = \Big(\sum_j
u_{ij}T_i^j\Big)^{-1} \sum_j du_{ij}T_i^j.
}{3.30}
Write
\eq{\Xi = \sum_{i,j}\lambda_{ij}\tau_{ij};\quad \lambda_{i,m_i-1}\neq 0
}{3.31}
where the nonvanishing condition comes from the requirement that the
form restricted to $C-D$ gives a trivialization along $\sD$ 
(see \eqref{2.30}). If we
write ($\mod T_i^{m_i}$)
\eq{\Big(u_{i0}+u_{i1}T_i+\ldots+u_{i,m_i-1}T_i^{m_i-1}\Big)^{-1} =
\nu_{i0}+ \nu_{i1}T_i+\ldots+\nu_{i,m_i-1}T_i^{m_i-1}
}{3.32} 
we get the table
\begin{align}\label{3.33}\tau_{i0} & = \nu_{i0}du_{i0} \\
\notag\tau_{i1} & = \nu_{i0}du_{i1}+\nu_{i1}du_{i0} \\
\notag\vdots \\
\notag\tau_{i,m_i-1} & = \nu_{i0}du_{i,m_i-1}+\ldots+\nu_{i,m_i-1}du_{i0}
\end{align}
Note if we give $u_{ij},\ du_{ij},\ \nu_{ij}$ all weight $j$, then
$\tau_{ij}$ will be homogeneous of weight $j$. 
Comparing \eqref{3.33} and \eqref{3.31}, it follows that if we expand
$\Xi|_B$ in terms of the $du_{ij}$, omitting $du_{10}$ and
$du_{1,m_1-1}$, we find for suitable $\alpha_{ij} \neq 0$
\ml{\Xi|_B = \sum_i g_{ij}du_{ij} = \sum_{i,j}
\Big[(\alpha_{i,m_i-1}u_{i0}^{-1}-\alpha_{1,m_1-1})du_{i,m_i-1}+ \\
u_{i0}^{-1}\sum_{p=0}^{m_i-2}
\Big(\alpha_{ip} \frac{u_{i,m_i-1-p}}{u_{i0}}+ 
\sum \text{ terms at least
quadratic in $\frac{u_{ik}}{u_{i0}}$} \Big)du_{ip}\Big]
}{3.34}
Looking at the weights, we see that for $p\le
m_i-2$
\ml{g_{ip} = \text{ nonzero multiple of
}\frac{u_{i,m_i-p-1}}{u_{i0}}+ \\
\text{ terms only involving $u_{ij},\
j<m_i-p-1$} 
}{3.35}
while
\eq{g_{i,m_i-1} = \alpha_{i,m_i-1}u_{i0}^{-1}-\alpha_{1,m_1-1}
}{3.36}
with neither $\alpha$ coefficient $0$. 

Now generators of $H^{m-2}_{DR}(B,(\sO,d+
\Xi))$ are of the form
$M\theta$ where $M$ is a monomial in the $u_{ij}, u_{i0}^{-1}$ and
$\theta$ is as in \eqref{3.29}. Relations are
\eq{\Big(\frac{\p}{\p u_{ij}}+g_{ij}\Big)(M)\theta = 0
}{3.37}
Because of \eqref{3.35} one can use these relations to eliminate
$u_{ij},\ j>0$ from $M$ by downward induction on $j$, starting from
$u_{i,m_i-1}$. We are left with the case $M=u_{2,0}^{n_2}\cdots
u_{r,0}^{n_r}$ with $n_i\in \Z$. In this case we can apply
\eqref{3.36}.  If $m_i\ge 2$, we get the relation
\eq{u_{2,0}^{n_2}\cdots u_{r,0}^{n_r}\theta \equiv
\frac{\alpha_{i,m_i-1}}{\alpha_{1,m_1-1}}u_{2,0}^{n_2}\cdots 
u_{i,0}^{n_i-1}\cdots u_{r,0}^{n_r} \theta  
}{3.38}
Using this, we can get $n_i=0$. If $m_i=1$ and $i\ge 2$ the relation
becomes
\eq{u_{2,0}^{n_2}\cdots u_{r,0}^{n_r}\theta \equiv
\frac{n_i+\alpha_{i,0}}{\alpha_{1,m_1-1}}u_{2,0}^{n_2}\cdots u_{i,0}^{n_i-1}\cdots
u_{r,0}^{n_r} \theta 
}{3.39}
If $\alpha_{i,0}$ is not a positive integer, we can arrange $n_i=0$. 

On the other hand, we claim that if $m_i=1$, then
$\alpha_{i0}$ is the residue of $i^*\nabla$ along $c_i$.
Indeed \eqref{2.29} shows that in this case,
$\tau_{i0}= d\log u_{i0}$, and $u_{i0}$ is then just the local
parameter in the point $c_i$ (see \eqref{2.28}).
Thus the quasiisomorphism 
\begin{multline}\label{3.40}
\{i^*\nabla_{/K}: i^*\sL \to i^*\sL \otimes
\Omega^1_{C/K}(\sD)\} \to \{j_*i^*\sL|(C-D) \\
\to j_*i^*\sL \otimes \Omega^1_{(C-D)/K}\} \end{multline}
forces $a_i$ not to lie in $\N -\{0\}$.
\end{proof} 

The following was left open in the proof of lemma \ref{lem3.8}:
\begin{lem}\label{lem3.10}Let $\Xi=\sum_{i,j}\lambda_{ij}\tau_{ij}$ be
as in lemma \ref{lem3.8}. Then $\Xi$ vanishes at a unique point $b\in
B$.  
\end{lem}
\begin{proof}We have seen in the proof of lemma \ref{lem3.8} that
$\Xi$ vanishes at a point in $B$. We must show it vanishes at at most
one point. Let $b=(\ldots,b_{ij},\ldots)\in B$ be a point. Write
$b=(\ldots,y_{ij},\ldots)$ with respect to the coordinates $\nu_{ij}$
\eqref{3.32}. Staring at \eqref{3.33}, the conditions that $\Xi|_b=0$
are seen to be (recall $\sum_i du_{i,m_i-1}=0=du_{10}$) for $i\ge 2$
\begin{gather}
\lambda_{i,m_i-1}y_{i1}+\lambda_{i,m_i-2}y_{10} =0  \label{3.41}\\
\lambda_{i,m_i-1}y_{i2}+\lambda_{i,m_i-2}y_{i1}+\lambda_{i,m_i-3}y_{i0}
 =0 \notag\\
\vdots \ \ \vdots \notag\\
\lambda_{i,m_i-1}y_{i,m_i-1}+\lambda_{i,m_i-2}y_{i,m_i-2}+\ldots +
\lambda_{i,0}y_{i0} = 0 \notag
\end{gather}
For $i=1$ one gets the same list but with the last line (coefficient
of $du_{i0}$) omitted. Finally, using $\nu_{10}=1$ and
$du_{1,m_1-1}=-\sum_{i\ge 2} du_{i,m_i-1}$ one gets
\eq{\lambda_{1,m_1-1}=y_{i0}\lambda_{i,m_i-1};\quad 2\le i\le r
}{3.42} 
Since $\lambda_{i,m_i-1}\neq 0$, equations \eqref{3.41} and
\eqref{3.42} admit a unique solution for the $y_{ij}$. Since we know
$\Xi$ vanishes at at least one point of $B$, this point must lie in
$B$. 
\end{proof}

Finally, we must calculate the Gau\ss-Manin connection on
$$H^{m-2}_{DR}(B,(\sO_B,d +\Xi))$$ 
Define $\Psi$ to be an absolute
invariant form lifting $\Xi$. By assumption our connection on $C-D$
comes from an absolute integrable connection which, by proposition
\ref{prop2.16}, comes from an absolute integrable connection on
$G$. Restricting this connection to $B$ gives our $\Psi$. 
\begin{lem}\label{lem3.11}With notation as above, there exists $F\in
\sO_B,\ \eta\in {\OK}$, and $a_i\in k,\ i\ge 2$, such that
\eq{\Psi = \sum_{i=2}^r a_i\frac{du_{i0}}{u_{i0}}+dF+\eta.
}{3.43}
If moreover $\nabla$ is integrable, then $\eta$ is closed. 
\end{lem}
\begin{proof}Since $\Xi$ is (relatively) closed on $B$, one can write
\eq{\Xi = \sum_{i\ge 2} a_i\frac{du_{i0}}{u_{i0}}+d_{/K}F;\quad a_i \in K.
}{3.44}
Lifting to an absolute form forces
\eq{\Psi = \sum_{i\ge 2} a_i\frac{du_{i0}}{u_{i0}}+dF + \sum_j
f_j\eta_j;\quad f_j\in \sO_B,\ \eta_j\in \OK.
}{3.45}
Here the $\eta_j$ are linearly independent in $\OK$.
Using $d\Psi=0$ modulo $\Omega^2_K \otimes \sO_B$
and taking residues along $u_{i0}=0$ yields $a_i\in
k\subset K$. Then computing $d\Psi \mod \sO_B\otimes\Omega_K^2$ yields
\eq{0=\sum_j d_{/K}f_j\otimes\eta_j\in \Omega^1_{B/K}\otimes\OK.
}{3.46}
It follows that $f_j\in K$, so $\eta:=\sum f_j\eta_j\in \OK$. Taking
$d$ again shows $\eta$ is closed if $\nabla$ is integrable.
\end{proof}

We now compute the Gau\ss-Manin connection. We have the diagram of
global sections
\eq{\minCDarrowwidth.5cm \begin{CD}@. \Omega^{m-2}_{B/k}
@>>\text{onto} >
\Omega^{m-2}_{B/K} \\
@. @VV d+\Psi V \\
\Omega^{m-2}_{B/K}\wedge \OK @>>\cong >
\frac{\Omega^{m-1}_{B/k}}{\Omega^{m-3}_{B/k}\cdot \Omega^2_K} @. 
\end{CD}
}{3.47}

The connection is determined by its value on $\theta$ \eqref{3.29}. To
calculate, one lifts $\theta$ to $\tilde\theta\in \Omega^{m-2}_{B/k}$
and then applies $d+\Psi$. But for $\tilde\theta$ one can choose the
form with the same expression \eqref{3.29}. This form is closed, so 
\eq{\nabla_{GM}(\theta) = \Psi\wedge\theta = (d_K(F)+\eta)\wedge\theta 
}{3.48}
Here we write $F=\sum_I a_Iu^I,\ a_i\in K$ and $d_K(F) := \sum
da_Iu^I$.  

Let $b\in B$ be the point corresponding to the trivialization of
$\Omega^1_{C/K}(\sD)$ given by the polar part of
the original relative connection. It really lies in $B$ since we
have assumed that ${\rm deg} \sL=0$. 
\begin{lem}\label{lem3.12} With notation as above, the 
Gau{\ss}-Manin connection on the rank $1$ $K$-vector space
$$H^{m-2}_{DR}(B,(\sO,d+\Xi))$$ 
described by \eqref{3.48} is isomorphic
to the connection on $K$ given by 
$$1\mapsto \Psi|_b + \frac{1}{2}d\log(\kappa)
$$
for a suitable $\kappa\in K^\times$. 
\end{lem}
\begin{proof}We have seen (lemma \ref{lem3.8}) that this point $b$ is
determined by the condition that $\Xi(b)=0\in \mathfrak
m_{B,b}/\mathfrak m_{B,b}^2$. Changing $\Psi$ by a closed form pulled back
from $K$ changes the Gau{\ss}-Manin connection and the connection at $b$
in the
same way, so we can assume $\eta=0$, i.e. $\Psi = \sum
a_i\frac{du_{i0}}{u_{i0}} +dF$. Write
\eq{g_{ij} = \begin{cases}\frac{\p
F}{\p u_{ij}} & j>0 \\ \frac{\p
F}{\p u_{i0}}+a_i/u_{i0} & j=0. \end{cases}
}{3.49}
Write $F = \sum_I a_I u^I$. Then 
\eq{\Psi = \sum g_{ij}du_{ij}+\sum_I
u^Ida_I;\quad \Psi\wedge\theta = \sum_I u^Ida_I\wedge\theta
}{3.50}
We have
\eq{g_{ij}(b) = 0,\ j<m_i-1;\quad g_{i,m_i-1}(b) = g_{k,m_k-1}(b);\
\text{all }i,k.
}{3.51}
Since $\sum_i du_{i,m_i-1}|_B=0$, we see from \eqref{3.51} that
\eq{\Psi|_{\{b\}} = \sum_Ib^Ida_I
}{3.52} 
Thus, it will suffice to relate $u^I\theta$ and $b^I\theta$ in
$H^{m-2}_{DR}$. Note that each monomial $u^I$ involves $u_{ij}$ for
only one value of $i$, and the weight of $u^I$ is $\le m_i-1$ (see the
discussion after \eqref{3.33}). 

Suppose first the weight of $u^I$ is strictly less than $m_i-1$. Let
$j$ be maximal such that $u_{ij}$ appears in $u^I$. From \eqref{3.35}
it follows that  
\eq{g_{i,m_i-1-j} = \alpha_{i,m_i-1-j}\frac{u_{ij}}{u_{i0}^2} + 
\text{terms involving only } u_{ik};\ k<j.
}{3.53} 
Here $\alpha_{i,m_i-1-j}\neq 0$. Define $u^L =
u^Iu_{i0}^2u_{ij}^{-1}$. Note the weight of $u^L$ is $<m_i-1-j$, so in
$H^{m-2}_{DR}$ we have (compare \eqref{3.37}) 
\ml{u^I\theta = (u^I - \alpha_{i,m_i-1-j}^{-1}(\frac{\p}{\p
u_{i,m_i-1-j}}+g_{i,m_i-1-j})u^L)\theta = \\
(u^I - \alpha_{i,m_i-1-j}^{-1}g_{i,m_i-1-j}u^L)\theta = \sigma_I\theta
}{3.54}
where $\sigma_I$ is a sum of terms
of weights $< |I|$ and terms of weight $|I|$ only involving
$u_{i0},\dotsc,u_{i,j-1}$. Note that $b^I= \sigma_I(b)$ because
$g_{i,m_i-1-j}(b) = 0$. In this way we reduce to the case
$u^I=u_{i0}^p$. Our assumption on the weight implies $m_i\ge 2$, so 
\eq{(\frac{\p}{\p u_{i,m_i-1}}+g_{i,m_i-1})u_{i0}^p =
g_{i,m_i-1}u_{i0}^p. 
}{3.55} 
Together with \eqref{3.35} and $g_{ij}(b)=0$, this enables us to
reduce to $p=0$.

Suppose now the weight of $I$ is
$m_i-1$. If $m_i\ge 2$ we can use the above argument, except in the
case $u^I = u_{ij}u_{i,m_i-1-j}u_{i0}^{-2}$. Here there are two subcases.
If $j\neq m_i-1-j$, the
$\alpha_{i,m_i-1-j}$ in 
\eqref{3.54} is $a_I$ in the expansion $F=\sum_I a_I u^I$, so
\eq{u^Ida_I\wedge\theta = (b^Ida_I+\frac{da_I}{a_I})\wedge\theta.
}{3.56}
This completes the proof in this case because the connections
$b^Ida_i$ and $b^Ida_I+d\log(a_I)$ are isomorphic. If, on the other hand,
$m_i$ is odd and $j=\frac{m_i-1}{2}$, the monomial $u^I =
u_{ij}^2u_{i0}^{-2}$ and $dF$ contains the term $2a_Iu_{ij}du_{ij}$.
Thus, from \eqref{3.53} we conclude $\alpha_{i,\frac{m_i-1}{2}}= 2a_I$.
The lefthand identity in \eqref{3.54} yields in this case
\eq{u^Ida_I\wedge\theta = (b^Ida_I + \frac{1}{2}d\log
a_I)\wedge\theta.
}{sp}
In the statement of the lemma, we take $\kappa$ to be the product of the
corresponding $a_I$. 

Suppose finally $m_i=1$. In this case $\frac{\p F}{\p
u_{i0}}=0$, so the corresponding $a_I=0$ and by \eqref{3.52} this term
contributes nothing to $\Psi|_{\{b\}}$. Similarly, by \eqref{3.50}
there is no contribution to $\Psi\wedge\theta$.  

\end{proof}

We give two interpretations of the $2$-torsion term
$\frac{1}{2}d\log(\kappa)$ occurring in the previous lemma. 

\begin{defn}\label{defins1} Let $\sigma$ be a closed $1$-form relative to
$K$ on
$$\Spec(K[[t_1,\dotsc,t_N]]).
$$
Assume $\sigma(0)=0\in \frak m/\frak m^2$.
Write $\sigma = dh$ with $h(0)=0$, so $h=h_2+h_3+\ldots$ with $h_i$
homogeneous of degree $i$. If $h_2$ is nondegenerate, we may define 
$\text{disc($\sigma)$}=
\text{discriminant}(h_2)\in K^{\times}/K^{\times 2}$.  This is
well-defined independent of the choice of parameters. 
\end{defn}
\begin{thm} The Gau{\ss}-Manin connection on
$H^{m-2}_{DR/K}(B,(\sO,d+\Xi))$ is isomorphic to 
$$(d+\Psi)|_{\{b\}} +
\frac{1}{2}d\log({\rm disc}(\Xi|_{\widehat{\sO}_{B,b}})).
$$
(In particular, the quadratic term in $h=\int
\Xi|_{\widehat{\sO}_{B,b}}$is non-degenerate.)
\end{thm}
\begin{proof}First we collect some facts about $\Xi =
\sum_{i,j}\lambda_{ij}\tau_{ij}= \sum g_{ij}du_{ij}$. We have $u_{1,0}=1$
and $\nu_{i,0} = u_{i,0}^{-1}$. It follows from \eqref{3.33} of the paper
that
$u_{1,m_1-1}$ does not appear in the expression for $\Xi$ and
$du_{1,m_1-1}$ only appears with constant coefficient. Restricting to
$B:u_{1,m_1-1} = -\sum_{i\ge 2} u_{i,m_i-1}$ thus has the effect of
surpressing the term in $du_{1,m_1-1}$ and changing the coefficients
$g_{i,m_i-1}$ by a constant for $i\ge 2$. Expressed in this way, it
follows that the coefficient of $du_{ij}$ in $\Xi|_B$ involves only
monomials in $u_{ip}$ for the same $i$. Giving $u_{ij}$ and $du_{ij}$
both weight $j$, the terms in $g_{ij}du_{ij}$ all have weight $\le
m_i-1$. It
follows from formulas \eqref{3.34}-\eqref{3.36} that, writing $U_{ij} =
u_{ij}-b_{ij}$ so $U_{ij}(b)=0$, we may write
$$\Xi|_B = \sum G_{ij}(U)dU_{ij}.
$$
Here $G_{ij}(0)=0$. giving $U_{ij}$ and $dU_{ij}$ weights $j$, all terms
with first index $i$ have weights $\le m_i-1$. All terms of the form
$$U_{ij}dU_{i,m_i-1-j},\quad 0\le j\le m_i-1,\ i\ge 2\quad(\text{resp.}\
i=1,\ 1\le j\le m_1-2)
$$
occur with nonzero coefficient. Notice that replacing $u_{ij}$ with
$U_{ij}+b_{ij}$ introduces monomials of lower degree, but these have
weight $< m_i-1$. 

It follows that $\text{disc}(\Xi|_{\widehat{\sO}_{B,b}})$ is the
determinant of a matrix 
$$M = \begin{pmatrix}M_1 & 0 &\hdots & \hdots & 0\\
0 & M_2 & 0 & \hdots & 0 \\
\vdots & \vdots & \vdots &\vdots & \vdots \\
0 & 0 & \hdots & \hdots & M_r
\end{pmatrix}
$$
where $M_i$ is symmetric, $m_i \times m_i$ (resp. $(m_1-2) \times
(m_1-2)$), and has the shape
$$\begin{pmatrix}\hdots &\hdots &\hdots &\hdots & \bullet \\
\hdots &\hdots &\hdots &\bullet & 0 \\
\hdots &\hdots &\bullet & 0 & 0 \\
\vdots &\vdots &\vdots &\vdots &\vdots \\
\bullet & 0 & \hdots & \hdots & 0
\end{pmatrix}
$$
with the entries $\bullet$ non-zero. 

Mod squares, $\det(M_i)$ is $1$ if $m_i$ is even, and is given by 
$$\frac{1}{2}\cdot \text{coefficient
of}(U_{i,\frac{m_i-1}{2}}dU_{i,\frac{m_i-1}{2}})
$$
if $m_i$ is odd. Writing 
\eq{\Psi|_B = \sum a_i \frac{du_{i0}}{u_{i0}} + dF 
}{ins1}
as just above \eqref{3.49} with $F=\sum a_I u^I$ we find 
$$\frac{1}{2} d\log\text{disc}(\Xi|_{\widehat{\sO}_{B,b}}) =
\frac{1}{2}\sum\frac{da_I}{a_I}= \frac{1}{2}d\log(\kappa).
$$

(The sum on the right is over all $I$ such that
$u^I = (u_{i,\frac{m_i-1}{2}})^2$.) 
\end{proof}

Another interpretation of the $2$-torsion is the following. As in
\eqref{2.27}  for $s_i$ a
local parameter at $c_i\in \sD$ and $t_i$ another copy of $s_i$ (so
$s_i-t_i$ is a local defining equation for the diagonal in $(C\times C$),
the pullback of
$u_{ip}$ to $K((s_i))$ is the coefficient $s_i^{-(p+1)}$ of $t_i^p$ in
$(s_i-t_i)^{-1}$. It follows from \eqref{3.30} that 
\begin{gather} \label{pbtau}
\tau_{ij} \ \text{pulls back
to} \
\frac{-ds_i}{s_i^{j+1}}.
\end{gather}

 Write the polar part of the connection at
$s_i=0$ in the form $(g_0+g_1s_i+\ldots)\frac{ds_i}{s_i^{m_i}}$. Since
$\Xi = \sum_{i,j}\lambda_{ij}\tau_{ij}$ pulls back to this connection
form, we get $g_0 = -\lambda_{i,m_i-1}$. On the other hand, again from
\eqref{3.30} the coefficient of
$u_{i,0}^{-2}u_{i,\frac{m_i-1}{2}}du_{i,\frac{m_i-1}{2}}$ in
$\Xi$ is $-\lambda_{i,m_i-1}$ if $m_i$ is odd. This coefficient is
the contribution to ${\rm disc}(\Xi|_{\sO_{B,b}})$ from the
point $c_i\in \sD$,
so we conclude

\begin{thm} Write the relative connection at $c_i\in \sD$ in the form
$(g_{i,0}+g_{i,1}s_i+\ldots)\frac{ds_i}{s_i^{m_i}}$. Then the Gau{\ss}-Manin
connection on
$H^{m-2}_{DR/K}(B,(\sO,d+\Xi))$ is isomorphic to 
$$(d+\Psi)|_{\{b\}} +
\sum_i \frac{m_i}{2}d\log(g_{i,0}(0)).
$$
\end{thm}

\begin{defn}\label{def3.16} With notation as above, write
$$\tau(L) = \sum_i \frac{m_i}{2}d\log(g_{i,0}(0)).
$$
\end{defn}

To summarize, we have proven

\begin{thm}\label{thm3.13}Let $C/K$ be a complete smooth curve of
genus $g$ over a field $K\supset k$. Let
$\nabla_{/K}:L \to L\otimes\OCK(\sD)$ 
be a connection, 
such that
$$\Big(L \to L\otimes\OCK(\sD)\Big) \to  
\Big(j_*L|(C-D) \to j_*L\otimes\OCK|(C-D)\Big) $$ 
is a quasiisomorphism. This implies that
the divisor $\sD$ is minimal such that $\nabla|_{C-D}$ extends with
values in $\OCK(\sD)$ (see section 4, \eqref{4.2}).
We also assume that $L$ has degree 0, and that
the connection on $L|_{C-D}$ lifts
to an integrable, absolute (i.e. $/k$) connection
$\tilde\nabla$. 
Then 
\eq{\nabla|_{\sD}:L|_\sD \to L\otimes\OCK(\sD)|_\sD
}{3.57} 
is an $\sO_\sD$-linear isomorphism and determines a trivialization of 
$$\OCK(\sD)|_\sD$$
Write $J_\sD$ for the generalized jacobian and
$J_{\sD,N}$ for the torseur of divisors of degree $N:=2g-2+\deg
\sD$. The above trivialization corresponds to a $K$-point $b\in
J_{\sD,N}$. Write $\pi_N : (C-D)^N \to J_{\sD,N}$ for the natural map,
and let $(L_N,\tilde\nabla_N)$ be the evident bundle and absolute
connection on $(C-D)^N$. Then there exists a unique invariant,
absolute connection $(\sL,\Phi)$ on $J_{\sD,N}$ such that
$\pi_N^*(\sL,\Phi)=(L_N,\tilde\nabla_N)$. Moreover, we have 
\eq{\Big(\det(H^*_{DR}(C-D,(L,\nabla)),\nabla_{GM}\Big)^{-1}
\cong (\sL,\Phi)|_{\{b\}} + \tau(L)
}{3.58}
where $\tau(L)$ is as in definition \ref{def3.16}.
\end{thm}
\begin{remark} \label{rmk3.14}
$2$-torsion also occurs in the determinant of de Rham cohomology for the
trivial connection \cite{BE}.  By virtue of the following lemma this can
only happen when the variety has even dimension. 
\end{remark}
\begin{lem} \label{lem3.15}
Let $X/K$ be a smooth projective variety of odd dimension $n=2m+1$ over a
function field in characteristic $0$. Then the Gau{\ss}-Manin determinant
$${\rm det}( H_{DR}(X/K), d)$$
is trivial in $\Omega^1_K/d\log K^*$.
\end{lem}
\begin{proof} The strong Lefschetz theorem identifies the determinant
connections on $H_{DR}^p$ and $H_{DR}^{2n-p}$ so we need only consider
the connection on $\det H_{DR}^{n}$. As well known, the Poincar\'e duality
morphism
$$\varphi: H^n_{DR}(X/K) \otimes H^n_{DR}(X/K)
\to H^{2n}_{DR}(X/K)= K$$ 
is compatible with the Gau{\ss}-Manin connection, which is
trivial on $H^{2n}_{DR}(X/K)= K$. 
On the other hand, it is alternating, thus its determinant
$${\rm det} (\varphi): {\rm det}(H^n_{DR}(X/K)) \otimes {\rm
det}(H^n_{DR}(X/K)) \to K$$
fulfills
$$ {\rm det} (\varphi)(e\otimes e)= p^2 \cdot 1$$
where $p\in K^*$ is the Pfaffian of the determinant
of $\varphi$, written in the basis $e$. Thus if 
$\nabla(e)=\alpha \otimes e$, one has 
$$ {\rm det}(\varphi)(\nabla(e\otimes e)) = 2 \alpha p^2 \cdot
1= 2p d(p) \cdot 1.$$
Thus $\alpha = d\log p$ and the determinant of the
Gau{\ss}-Manin connection is trivial.
\end{proof}
\section{Product and Trace}

In this section, we introduce a product which is reminiscent of
Deligne's product explained in \cite{Dsm}.

We keep the notations of  sections 2 and 3
for $C/K,\ (L, \nabla),\ j:U=C-D \to C$, $\sD=\sum
m_ic_i$,  and $\sD' = \sD - D$. Further, 
\begin{gather}\label{4.1}
\nabla: L \to L\otimes \Omega^1_{U}
\end{gather}
is an absolute connection with vertical curvature $\nabla^2(L)\subset
L\otimes\Omega^2_K\otimes K(X)$.  Let $\nabla_{/K}: \sL \to \sL \otimes
\Omega^1_{C/K}(\sD)$ be an extension of $(L, \nabla_{/K})$ such that  
\begin{gather}\label{4.2}
\{\sL \to \sL
\otimes \Omega^1_{C/K}(\sD)\} \to
\{j_*L \to j_*L\otimes \Omega^1_{U/K}\}
\end{gather}
is a quasiisomorphism. We assume $\nabla_{/K}$ has a pole at every $c_i\in
D$. Note this implies that $\nabla_{/K}$ does not factor through 
$\Omega^1_{C/K}(\sD-c_i)$ for any $i$. Indeed, by assumption, the complex
\eq{j_*L/\sL \to (j_*L/\sL)\otimes\OCK(\sD)
}{4.3}
is acyclic. Take $e$ a local basis of $\sL$ at $c_i$ and $z$ a local
parameter, and suppose the connection can be written locally as 
$\nabla_{/K}e =
a(z)dz/z^{m-1}e$ with $m=m_i$. Then $\nabla_{/K}(z^{-1}e) =
(a(z)dz/z^{m}-z^{-2}dz)e$. The assumption that $\nabla_{/K}$
does have a pole at 
$c_i$ implies   that $m\ge 2$, 
so $z^{-1}e$ would represent a nontrivial element in
$H^0$ of the complex \eqref{4.3}, a contradiction. 

By lemma \ref{lem3.1} we know that the verticality condition
implies that the absolute connection extends as
\begin{gather}\label{4.3bis}
\nabla: \sL \to \sL\otimes
\Omega^1_{C}\ldr(\sD').
\end{gather}
{From} now on, we fix such a $(\sL, \nabla)$.

As we have seen, the map $\sL|_\sD \to \sL\otimes\OCK(\sD)|_\sD$ is function
linear. Since the connection does not factor through lower order poles, this
gives a trivialization (denoted $\text{triv}(\nabla)$) of $\OCK(\sD)|_\sD$. We
have   
\begin{gather}\label{4.4}
(c_1(\Omega^1_{C/K}(\sD)), {\rm triv \ }(\nabla)) \in 
\H^1(C, \sO^*_C \to \sO^*_{\sD}) \\
(\sL, \nabla) \in \H^1(C, \sO^*_C \to \Omega^1_C\ldr(\sD')).\notag
\end{gather}
The aim of this section is to define a product
\begin{gather}\label{4.5}
\cup: \H^1(C, \sO^*_C \to \sO^*_{\sD}) \times
\H^1(C, \sO^*_C \to \Omega^1_C\ldr(\sD')) \notag \\
\to \H^2(C, \sK_2 \to \Omega^2_C)
\end{gather}
Here $\sK_2$ is the Milnor sheaf associated to $K_2$, and the map $\sK_2 \to
\Omega^2_C$ is the $d\log$ map $\{a,b\} \mapsto
\frac{da}{a}\wedge\frac{db}{b}$. For a more detailed study of characteristic
classes for connections defined in the hypercohomology of such complexes, the
reader is referred to \cite{E}. In addition, we will define a trace
\begin{gather}\label{4.6}
{\rm Tr \ }:\H^2(C, \sK_2 \to \Omega^2_C) \to \Omega^1_K/d\log K^*
\end{gather}
We write 
\eq{A\cdot B := \text{Tr}(A\cup B)
}{4.6bis} 
so for example
\begin{multline}\label{4.7}
(c_1(\Omega^1_{C/K}(\sD)), {\rm triv \ }(\nabla))\cdot (\sL, \nabla)
:=  \\
{\rm Tr \ }((c_1(\Omega^1_{C/K}(\sD)), {\rm triv \ }(\nabla))\cup
(\sL, \nabla)) 
\end{multline}

Let 
$\sO^*_{C, \sD}= {\rm Ker \ }  (\sO^*_C \to \sO^*_{\sD})$.
Then 
\begin{lem}\label{lem4.1}
$d \log \sO^*_{C, \sD} \wedge \Omega^1_C\ldr(\sD') \subset 
\Omega^2_C$
\end{lem}
\begin{proof}
Since $\sO^*_{C, \sD}\subset 1+ \sI_{\sD}$, where
$\sI_{\sD}$ is the ideal sheaf of $\sD$, 
$$d\log 
\sO^*_{C, \sD} \subset \sO_C d \sI_{\sD} \subset \Omega^1_C(*D).$$ 
Also one has
$d \sI_{\sD} \subset \sI_{\sD}\otimes_{\sO_C}\Omega^1_C\ldr$.
Thus $$ 
d \log \sO^*_{C, \sD} \wedge \Omega^1_C\ldr(\sD')
\subset \sI_D \otimes_{\sO_C}\Omega^2_C\ldr \subset \Omega^2_C.$$
\end{proof}
We define $\cup$ by
\begin{gather}\label{4.8}
\sO^*_{C, \sD} \cup \sO^*_C \to \sK_2  \\
(\lambda, c) \mapsto \{\lambda, c\} \notag \\
\sO^*_{C, \sD} \cup \Omega^1_C\ldr(\sD') \to \Omega^2_C \notag \\
(\lambda, \omega) \mapsto d \log \lambda \wedge \omega.\notag
\end{gather}
Concretely, we can write the product in terms of Cech cocyles. Here
$\sC^i$ refers to Cech cochains, $\delta$ is the Cech coboundary, and $d$ is a
boundary in the complex:
\begin{gather}
(\lambda_{ij}, \mu_i) \in (\sC^1(\sO^*_C) \times
\sC^0(\sO^*_{\sD}))_{d-\delta} \notag \\
(c_{ij}, \omega_i) \in (\sC^1(\sO^*_C) \times
\sC^0(\Omega^1_C\ldr(\sD'))_{d-\delta}\notag
\end{gather}
one has
\begin{gather}\label{4.9}
(\lambda, \mu) \cup (c, \omega) =
(\{\lambda_{ij}, c_{jk}\}, d\log \lambda_{ij} \wedge \omega_j, 
-d \log \tilde{\mu_i} \wedge \omega_i )  \\
\in (\sC^2(\sK_2) \times \sC^1(\Omega^2_C\ldr(\sD')) \times
\sC^0(\Omega^2_C\ldr(\sD')/\Omega^2_C))_{d+\delta}\notag 
\end{gather}
where $\tilde{\mu_i} \in \sC^0(\sO^*_C)$ is a local lifting of
$\mu_i$. Note we have replaced the complex $\sK_2 \to \Omega^2_C$ with the
quasiisomorphic complex 
$$\sK_2 \to \Omega^2_C\ldr(\sD') \to
\Omega^2_C\ldr(\sD')/\Omega^2_C.
$$
\begin{prop}\label{prop4.2}
The product $\cup$ extends to
\begin{gather}
\cup: \H^1(C, \sO^*_C \to \sO^*_{\sD}) \times
\H^1(C, j_*\sO^*_U \to \Omega^1_C\ldr(\sD')) \notag \\
\to \H^2(C, \sK_2 \to \Omega^2_C).\notag
\end{gather}
\end{prop}
\begin{proof}
The map 
$$\H^1(C, \sO^*_C \to \Omega^1_C\ldr(\sD'))
\to \H^1(C, j_*\sO^*_U \to \Omega^1_C\ldr(\sD'))$$ 
is surjective,
and its kernel is the $\Z$-module generated by $(\sO(D_i),
d_i)$ where $d_i$ is the connection with logarithmic poles along
$D_i$ with residue -1. 
Let $z_1$ be a local coordinate around $c_1$. 
Let $U_i$ be a Cech covering of $C$,
with  $c_1 \in U_1 \subset V_1$, and $c_1 \notin U_i, i \neq 1$. Assume $c_1$
is the only zero or pole of $z_1$ on $U_1$.  Let 
$$(\lambda, \mu) \in \H^1(C,\sO^*_C\to\sO^*_{\sD})
$$ 
be a Cech representative of a class in $ \H^1(C, \sO^*_C \to
\sO^*_{\sD})$. Then $(c_{ij}, \omega_i)$ with
$c_{1j}= z_1^{-1}$, $c_{ij}= 1$ for $i \neq 1$, $\omega_1=-d
\log z_1$, $\omega_i =0$ for $i\neq 1$ is a Cech representative
of  $(\sO(D_1),d_1)$. 
Thus considering $Z \in \sC^0(\sO_C[z_1^{-1}]^*)$ with $Z_1=z_1$ and
$Z_i=1$ for $i\neq 1$,
the cocyle of \eqref{4.6} is just the coboundary
$$(d-\delta)(\{\lambda_{ij}, Z_j \}, d\log \tilde{\mu_i} \wedge d \log
Z_i) \in (d-\delta)(\sC^1(\sK_2) \times \sC^0(\Omega^2_C\ldr(\sD')).$$
(Note $Z_j$ is invertible on $U_{ij}$ for $i\neq j$ so the $K_2$-cochain is
defined.) 
\end{proof}
Now we define the trace. We have (with standard $K$-theoretic
notation, \cite{B}) 
\begin{gather}\label{4.10}
H^2(C, \sK_2)=0 \\
{\rm Nm \ } : H^1(C, \sK_2) = \{\oplus _{x\in C^{(1)}} \kappa(x)^*\}/
{\rm Tame}(K_2(K(C))) \to K^* \notag\\
\sum_x \varphi_x \mapsto \Pi_x {\rm Nm \ }(\varphi_x) \notag 
\end{gather} 
and of course $H^1(C, \Omega^2_C)= \Omega^1_K \otimes H^1(C,
\Omega^1_{C/K})= \Omega^1_K$.
This defines
\begin{gather}\label{4.11}
{\rm Tr \ }:\H^2(C, \sK_2 \to \Omega^2_C)= H^1(C, \Omega^2_C)/H^1(C,
\sK_2)
\to \Omega^1_K/d \log K^* .
\end{gather} 
\begin{lem}\label{lem4.3}
The trace
\begin{gather}
{\rm Tr \ }:\H^2(C, \sK_2 \to \Omega^2_C)=
\H^2(C, \sK_2 \to \Omega^2_C\ldr(\sD') \to
\Omega^2_C\ldr(\sD')/\Omega^2_C) \notag \\
\to \Omega^1_K/d\log K^* \notag
\end{gather}
factors through
\begin{gather}\label{4.13a}
\H^2\Big(C, \sK_2 \to \Omega^1_K \otimes \Omega^1_{C/K}(\sD) \to
\Omega^1_K \otimes (\Omega^1_{C/K}(\sD)/\Omega^1_{C/K})\Big) \\
\cong \Omega^1_K/d \log K^* \notag \\ 
\cong \Omega^1_K \otimes_K H^0(\sD, \omega_{\sD/K})/
\H^1(C, \sK_2 \to \Omega^1_K \otimes \Omega^1_{C/K}(\sD)) \notag
\end{gather}
where $\omega_{\sD/K}$ is the relative dualizing sheaf of the
scheme $\sD$, containing $K \cong \omega_{D/K}$.
\end{lem}
\begin{proof}Note that
$$\OK\otimes_K \OCK(\sD) \cong
\Omega^2_C\ldr(\sD')/(\Omega^2_K\otimes\sO_C(\sD') )
$$
so the complex in \eqref{4.13a} is indeed a quotient. From the diagram
$$\minCDarrowwidth.5cm\begin{CD}@. \sK_2 @= \sK_2 @. @. \\
@. @VVV @VVV @. \\
0 @>>> \OK\otimes\OCK @>>> \OK\otimes\OCK(\sD) @>>> \OK\otimes(\OCK(\sD)/\OCK)
@>>> 0 
\end{CD}
$$
one deduces that the left hand side of \eqref{4.13a} is isomorphic to
$$\H^2(C,\sK_2 \to \OK\otimes\OCK)\cong \text{coker}(H^1(C,\sK_2) \to \OK
\otimes H^1(C,\OCK)). 
$$
The right hand side here is identified under the norm with
$\OK/d \log K^*$, which proves the second equality. The third
one comes from the map $$\Omega^1_K\otimes \omega_{\sD/K}[-2]
\to \{\sK_2 \to \Omega^1_K \otimes \Omega^1_{C/K}(\sD) \to
\Omega^1_K \otimes \omega_{\sD/K}\}$$ and the vanishing of
$\H^2(\sK_2 \Omega^1_K \otimes \Omega^1_C(\sD))$. Note that this
cumbersome way of writing this cohomology allows to write local
contribution of a class in this cohomology group.
\end{proof}

The first main result of this section is the following
\begin{thm}\label{thm4.4}
Let $(\sL, \nabla)$ and $(\sL', \nabla')$ be two extensions of
the vertical connection 
$(L, \nabla)$ on $U$ as above satisfying the quasiisomorphism
condition \eqref{4.2}. Then, with notation as in \eqref{4.8},
$$((c_1(\Omega^1_{C/K}), {\rm triv \ } \nabla)\cdot (\sL,
\nabla))= ((c_1(\Omega^1_{C/K}), {\rm triv \ } \nabla') \cdot (\sL',
\nabla')).$$
\end{thm}
\begin{proof}

The quasiisomorphism condition is local about each point of $D$, so we
may assume our line bundles are $\sL(\nu c)\subset \sL$ for some $\nu
< 0$ and $c\in D$.  

Choose local coordinates $z_i$ near $c_i$ and a Cech covering
$U_i$ of $C$ such that $c_i \in U_i$, $z_i \in \sO^*(U_i -c_i)$,
$c_i \notin U_j$ for
$i\neq j$. Let us denote by 
$$(c_1(\Omega^1_{C/K}), z_i) \in \H^1(C, \sO^*_C \to
\sO^*_{\sD})$$ the class defined by the local trivialization
$$\frac{dz_i}{z_i^{m_i}}: \sO_{m_ic_i} \to \Omega^1_{C/K}(\sD)
\otimes \sO_{m_ic_i}.$$
Let $(\sL, \nabla) = (c_{ij}, \omega_i)$. Then 
$\omega_i= a_i \frac{dz_i}{z_i^{m_i}} +
\frac{b_i}{z_i^{m_i-1}},$ with $a_i \in \sO_C$ such that
$a_i|_{m_ic_i} \in \sO^*_{m_ic_i}$ and
$b_i \in \Omega^1_K \otimes \sO_C.$ Suppose $c=c_i$. We drop the index
$i$ for convenience. One has
\begin{gather}\label{bigform}
((c_1(\Omega^1_{C/K}), {\rm triv \ } \nabla)\cdot (\sL,
\nabla))= \\
(c_1(\Omega^1_{C/K}), z)\cdot (\sL,
\nabla) + \Big(0, 0,- d\log (a) \wedge (a \frac{dz}{z^{m}} +
\frac{b}{z^{m-1}})\Big). \notag
\end{gather}
where the last term is a cocycle as in \eqref{4.10} or the quotient
complex \eqref{4.13}. For $(\sL(\nu c),\nabla(\nu c))$ one replaces 
$a$ by $a-\nu z^{m-1}$, 
leaving $b$ and $m$ unchanged. 

By proposition \ref{prop4.2}, one has
\begin{gather}
(c_1(\Omega^1_{C/K}), z)\cdot (\sL,
\nabla)= (c_1(\Omega^1_{C/K}), z)\cdot (\sL(\nu c),
\nabla(\nu c)) \notag
\end{gather}
Thus
\begin{gather}
((c_1(\Omega^1_{C/K}), {\rm triv \ } \nabla)\cdot (\sL,
\nabla)) -((c_1(\Omega^1_{C/K}), {\rm triv \ } \nabla') \cdot (\sL',
\nabla')) \notag \\
= \Big(0,0, (d(a- \nu z^{m-1})  - d(a))\wedge \frac{dz}{z^m} + d\log
(\frac{a - \nu z^{m-1}}{a}) \wedge \frac{b}{z^{m-1}}\Big)
\notag \\
= \Big(0,0,\frac{ \nu (m-1) d\log z \wedge \frac{b}{a}}{(1-\nu 
\frac{z^{m-1}}{a})}\Big).\notag 
\end{gather}
The nontrivial part in the last expression is computed in
$$H^0(\Omega^2_C\langle c\rangle/\Omega^2_C)\cong \Omega^1_c.
$$
Computing using the residue at $c$ we find the above difference is
$$\Big(0,0,\nu (m-1)\frac{b(c)}{a(c)}\Big)
$$
The verticality condition for the curvature reads
$$
$$
$$ da \wedge \frac{dz}{z^m} + \frac{db}{z^{m-1}}
-(m-1)\frac{dz}{z^m} \wedge b =0\in \Omega^1_K \otimes
\Omega^1_{C/K}(\sD). 
$$  
In particular, $(da - (m-1) b)_{|c} = 0$.  
The difference of the two products
is therefore $(0,0, d\log a^\nu)$, which vanishes in $\Omega^1_K /d\log K^*$.
\end{proof}   
\begin{remark}A version of the formula \eqref{bigform} in higher rank plays a
central role in section 5.
\end{remark}

Suppose now
\begin{gather}\label{4.12}
m_i=1 {\rm \ for \ all \ } i.
\end{gather}
In this case, the class $
(c_1(\Omega^1_{C/K}(D), z_i) \in \H^1(C, \sO^*_C \to \sO^*_{\sD})$
as defined in the proof of theorem \ref{thm4.4} does not in fact depend
on the choice of the local coordinate $z_i$. Indeed, the
trivialization $$\sO_D \to (\Omega^1_{C/K}(D)/\Omega^1_{C/K}=
\sO_D), 1 \mapsto \frac{dz_i}{z_i}$$ is just the canonical
identification given by the residue along $c_i$. In other words,
the class $(c_1(\Omega^1_{C/K}(D), z_i)$ is what is denoted by 
$(c_1(\Omega^1_{C/K}(D), {\rm res }_D)$ 
in \cite{BE}, and appears on the right hand side of the
Riemann-Roch formula.
The second main result of this section is
\begin{thm}\label{thm4.5}
Let $(\sL, \nabla)$ be as above, with $m_i=1$ for
all $i$. Then 
\begin{multline*}{\rm det \ } \Big(H^*_{DR}(U, L)),\ \text{Gau\ss\ - Manin
connection}\Big) = \\
- c_1\Big(\Omega^1_{C/K}(D), {\rm
triv }\nabla \Big)\cdot (\sL, \nabla).
\end{multline*}
\end{thm}
\begin{proof}
Given the main result of \cite{BE}, and lemma \ref{lem3.15},
the theorem is of course
equivalent to
\begin{gather}\label{4.13}
c_1(\Omega^1_{C/K}(D), {\rm
res}_D) \cdot (\sL, \nabla) = c_1(\Omega^1_{C/K}(D), {\rm triv
\ } \nabla) \cdot (\sL, \nabla).
\end{gather}
Keeping the same notations as in the proof of theorem
\ref{thm4.4}, one has
\begin{gather}
c_1(\Omega^1_{C/K}(D), {\rm
res}_D) - c_1(\Omega^1_{C/K}(D), {\rm
triv  \ } \nabla)= (0, a_i) \notag
\end{gather}
 and thus 
\begin{gather}
(c_1(\Omega^1_{C/K}(D), {\rm res}_D) - c_1(\Omega^1_{C/K}(D), {\rm
triv  \ } \nabla)) \cdot (\sL, \nabla) = \notag \\
(0,0, -d (a_i) \wedge d\log z_i - d\log a_i \wedge b_i) = \notag
\\ (0,0, -d(a_i) \wedge d\log z_i). \notag
\end{gather}  
This lies in $\Omega^1_K \otimes \omega_{D/K}=\Omega^1_K$ and by
lemma \ref{lem4.3}, its trace factors through $\Omega^1_K
\otimes H^1(C, \Omega^1_{C/K})$. But the image of $\gamma=\sum_i  a_i
d\log z_i \in \omega_{D/K}$ in $H^1(C, \Omega^1_{C/K})$ is the 
relative Atiyah class 
${\rm at }_{/K}(\sL)$ (\cite{EV}, appendix B), thus the image of
$d\gamma= \sum_i  d(a_i)\wedge
d\log z_i \in \Omega^1_K \otimes \omega_{D/K}$
in $\Omega^1_K \otimes H^1(C, \Omega^1_{C/K})$ 
is $d({\rm at }(\sL))$, where ${\rm at}(\sL) \in H^1(C,
\Omega^1_C)$ is the absolute Atiyah class of $\sL$. Indeed, 
$d: H^1(C, \Omega^1_C) \to \Omega^1_K \otimes H^1(C,
\Omega^1_{C/K})$ factors through $H^1(C, \Omega^1_{C/K})$ by
Hodge theory. On the other hand, if $c_{ij} \in \sC^1(\sO^*_C)$
is a cocyle for $\sL$, then $d\log c_{ij} \in \sC^1(\Omega^1_C)$
is a cocyle for ${\rm at}(\sL)$, and consequently, $d({\rm at }(\sL))=0$.
\end{proof} 

We want to explain briefly a fundamental compatibility satisfied by the pairing
\eqref{4.6bis}. Let $b = \sum b_i$ be a $0$-cycle on $C$ with support disjoint
from $\sD$, and let $\nabla:\sL\to
\sL\otimes\Omega^1_{C}<D>(\sD')$ be an absolute,
integrable connection. We can interpret $\nabla|_{b_i} \in
\Omega^1_{K(b_i)}/d\log K(b_i)^*$. 
\begin{prop}\label{prop4.6} With notation as above, let $[b]\in
\H^1(C,\sO_C^* \to \sO_\sD^*)$ be the class of the $0$-cycle $b$. Then
$$[b]\cdot (\sL,\nabla) = \sum_i
{\rm Tr}_{K(b_i)/K}(\nabla|_{b_i})\in \OK/d\log K^*. 
$$
Let $(\sL_0, \nabla_0)$ be the invariant connection on $J_\D$
which pulls back to $(\sL, \nabla)$ via the cycle map $i: (C-D)
\to J_\sD$ (proposition \ref{prop2.16}). Then 
$$[b]\cdot (\sL,\nabla)= {\rm Tr}_{i_0(b)/} (\sL_0, \nabla_0)|_{i_0(b)}.
$$
\end{prop}
\begin{proof} 
One reduces easily to the case $b$ is a single $K$-point. Let
$U_2$ be a Zariski-open set containing $D$ and $b$. Shrink $U_2$
if necessary so
there exists $z\in H^0(\sO_{U_2})$ with $z|_\sD = 1$ and $(z)=b$.  Let
$U_1=C-\{b\}-\sD$ so $C=U_1\cup U_2$.  Shrinking the $U_i$ if necessary, we can
assume $\sL|_{U_i}\cong \sO_{U_i}$, so $(\sL,\nabla)$ is represented by some
cocycle $(\mu_{12},\omega_1,\omega_2)$. Then
$$\nabla|_b = \omega_2|_b \in \OK/d\log K^*
$$
On the other hand, by the definition \eqref{4.10},
$[b]\cdot (\sL,\nabla)$ is represented by the image of the cocycle
$$d\log(z)\wedge \omega_2|_{U_{12}}
\in H^1(C,\Omega^2_C)
\to
\OK\otimes_K H^1(C,\OCK) \cong \OK.
$$
Write $\omega_2 = \omega_2(b)+z\eta_2$ with $\eta_2$ regular on $U_2$. Since
$d\log(z)\wedge z\eta_2$ extends to $U_2$, it is homologous to zero, so 
$$[b]\cdot (\sL,\nabla) = \omega_2(b)[b]\in \OK\otimes\H^1(C,\OCK) \mapsto
\omega_2(b)\in \OK.
$$
Now, since one obviously has
$${\rm Tr} (\nabla|_{b_i}) = {\rm Tr} (\nabla_0)|_{i_0(b_i)}
$$
and the translation $\nabla_0$ is invariant, the second
equality is a direct interpretation of the first one.
\end{proof}

Now we can formulate and prove a variant of theorem \ref{thm3.13}.
\begin{thm}\label{thm4.7}
Let $(C/K, U/K, (L, \nabla),\sD)$ be as in \eqref{4.1},
\eqref{4.2}, \eqref{4.3}. Then
\begin{gather}\label{4.14}
{\rm det \ }(H^*_{DR}(U, L))= - (c_1(\Omega^1_{C/K}(\sD), {\rm
triv \ }(\nabla))\cdot (\sL, \nabla) \in \Omega^1_K/d\log K^*
\end{gather}
modulo torsion (see remark \ref{rmk3.14}), 
where $\nabla_{/K}: \sL \to \sL \otimes \Omega^1_{C/K}(\sD)$ is
any extension of $(L, \nabla_{/K})$ having poles along all points of
$\sD$ such that
\begin{gather}
\{\sL \to \sL \otimes \Omega^1_{C/K}(\sD) \} \to \{j_*L \to
j_*(L\otimes \Omega^1_{U/K})\}\notag
\end{gather}
is a quasiisomorphism.
\end{thm}
\begin{proof}
If all $m_i=1$ this is just theorem \ref{thm4.5}. So we assume
that $m_1 \geq 2$ in the sequel. Then as in the proof of theorem
\ref{thm4.4}, replacing $\sL$ by $\sL(-c_1)$ changes $a_1$ to
$(a_1 +z_1^{m_1-1})$ and keeps the rest unchanged. Thus
the quotient complex 
$$\sL/\sL(-c_1) \to \sL/\sL(-c_1) \otimes \Omega^1_{C/K}(\sD)$$
is $\sO_{c_1}$-linear and the map is the multiplication by $a_1
\in \sO^*_{c_1}$. In particular, $\sL(\nu c_1)$ fulfills \eqref{4.2} for
all $\nu \in \Z$ and taking $\nu=- {\rm deg \ } \sL$, we may
assume by theorem \ref{thm4.4} that ${\rm deg \ } \sL= 0$. 
If $H^0_{DR}(U,L) \neq 0$, then there is a meromorphic section $\varphi$ of
$\sL$ verifying the flatness condition
$$d\varphi + \omega \varphi=0.$$
This implies in particular that $\omega$ has at most logarithmic
poles along $D$, which contradicts the condition $m_1 \geq 2$. 
On the other hand, $H^2_{DR}(U, L) =0$ for dimension reasons,
thus we can apply theorem
\ref{thm3.13} together with 
proposition \ref{prop4.6}  to obtain the result, after we
have replaced $(\sL, \nabla)$ by $(\sL, \nabla) \otimes f^*\Big((\sL,
\nabla)|_{c_0}^{-1}\Big)$ to trivialize the 
connection at $c_0$ and applied the projection formula
to this tensor connection.
\end{proof}
We finish this section with an example. Let $U=\A^1_K$ be the
affine line over ${\rm Spec}K$, with parameter $t$,
and $\nabla$ be a connection on the trivial bundle. Then 
up to a twist by a form of the base, $\nabla$ has equation
$A= df$, where $f= \sum_{i=1}^{m-1}a_i t^i$, $a_i \in K, a_{m-1}
\neq 0$. Write $df= d_Kf + f'dt$ with $d_kf= \sum_{i=1}^{m-2}
da_i t^i$ and $f'= \sum_{i=1}^{m-2} ia_i t^{i-1}$. Let $b_i,
i=1,\ldots, (m-2)$ be the zeroes of $f'$ (defined over $K$ after
some finite field extension), and let $N_\ell(\underline{b})=
\sum_{i=1}^{m-2} b_i^\ell$ be the Newton classes of the zeroes
of $f'$, which of course are expressable in the $a_i$ already on
${\rm Spec}K$. Then the main theorem says
$${\rm det}(GM)^{-1}= \sum_{i=1}^{m-2} da_i N_i(\underline{b}).$$

\section{A Formula in Higher Rank}

In this section, we want to define a sort of non-commutative
product of a higher rank connection with the Chern class of 
the dualizing sheaf of $C$ with poles. The notations are as in
the whole article: $C$ is a curve defined over a function field
$K$ over an algebraically closed field $k$ of characteristic 0,
$U$ is an open set such that $D=X-U=\sum_i c_i$ consists of $K-$rational
points. Let $(\sE, \nabla)$ be a rank $r$-connection on $U$ with vertical
curvature \eqref{2.60}. Let $m_i$ be the multiplicity of the relative
connection at the point $c_i$, that is, the minimal multiplicity
such $\nabla $ factors
$$\nabla_{/K}: E \to E\otimes \Omega^1_{C/K}(\sum_i  m_ic_i).$$
Lemma \ref{lem3.1} no longer holds
true in the higher rank case. 
We say that the poles of the global connection {\it behave
well} if
\begin{gather} \label{ass5.1}
\nabla: E \to E \otimes \Omega^1_C<D>(\sD')
\end{gather}
where $\sD = \sum_i m_i c_i$ and $\sD'= \sD - D$.

Let $s=\{s_i\}$ be a trivialization of $\omega_{C/K}(\sD)|_\sD
\cong \omega_{\sD/K}$. That
is, $s_i\in \omega_{m_ic_i}$ and the map $1\mapsto s_i$ is an isomorphism
$\sO_\sD \cong \omega_{m_ic_i}$. For example, if $z_i$ is a local parameter,
one can take $s_i = \frac{dz_i}{z_i^{m_i}}$. We will abuse notation and write
$s_i$ also for a lifting of the trivialization to a local section of
$\Omega^1_C< D>(\sD')$. The local matrix of the connection has the shape
\begin{gather}
A_i = g_i s_i + \frac{\eta_i}{z_i^{m_i-1}}
\end{gather}
where $g_i$ and $\eta_i$ are $r\times r$ matrices with
coefficients in $\sO_C$ and $\Omega^1_K \otimes \sO_C$ 
respectively. Note that the matrix of functions $g_i$ depends only on the
lifting of $s_i$ to a section of $\omega_{C/K}(\sD)$. (Indeed, the relative
connection has matrix $gs_i$.)

We assume
\begin{gather} 
{\rm Image}(g_i) \in M(r \times r, \sO_{\sD}) \notag \\
\label{ass5.3}
{\rm lies  \  in \ } GL(r,  \sO_{\sD})
\end{gather}
Under this assumption, we define
\begin{gather}\label{def5.4}
\{c_1(\omega(\sD)), \nabla\}:=  \\
 c_1(\omega(\sD),
s_i) \cdot {\rm det}( \nabla) - \sum_i{\rm res \
Tr}(dg_i g_i^{-1}A_i) \notag \\
\in \Omega^1_K/d\log K^* \notag
\end{gather}
\begin{conj} \label{conj} Assuming \eqref{ass5.1} and
\eqref{ass5.3}, we have
$${\rm det}H^*_{DR}(U, (\sE, \nabla))^{-1} = 
\{c_1(\omega(\sD)), \nabla\} \in (\Omega^1_K/d\log K^*) \otimes_{\Z} \Q.
$$
\end{conj}
We discuss the assumption \eqref{ass5.3} (see proposition
\ref{prop5.6}) at the end of this section. 
The assumption \eqref{ass5.1} on the poles behaving
well is not very well understood. It reflects a sort of
stability in all possible directions for
the poles of the global connection.

First, we justify   
the conjecture by establishing some rather surprising invariance
properties for $\{c_1(\omega(\sD)), \nabla\}$.
\begin{lem}\label{lem5.2}Fix an index $i$ and write the connection matrix
locally in the form
$$A = g\frac{dz}{z^m}+\frac{\eta}{z^{m-1}}
$$
where $g$ is an invertible matrix of functions and $\eta$ is a matrix with
entries from $\Omega^1_K\otimes\sO_C$. Then \begin{enumerate}\item $\text{res
Tr}(dgg^{-1}A)= \text{res Tr}(dgg^{-1}\frac{\eta}{z^{m-1}})$.
\item $[\eta,g]z^{1-m}$ has no pole at the point $z=0$.
\end{enumerate}
\end{lem}
\begin{proof}The assumption that the curvature is vertical implies
\begin{equation} \label{eq5.5}
dA = dg\frac{dz}{z^m}+d(\frac{\eta}{z^{m-1}}) \equiv A^2 =
[\eta,g]\frac{dz}{z^{2m-1}} \mod\ \Omega^2_K\otimes\sO_C[z^{-1}]
\end{equation}
Multiplying through by $z^m$ and contracting against $\frac{\partial}{\partial
z}$ we deduce 2. 

For 1, we must show $\text{res Tr}(dg\frac{dz}{z^m}) = 0$. From \ref{eq5.5},
using $\text{Tr}[g,\eta] = 0$ we reduce to showing $\text{res
Tr}d(\frac{\eta}{z^{m-1}}) = 0$. 
Since $\eta$ has entries  $\Omega^1_K$, one has
$$\text{res
Tr}d(\frac{\eta}{z^{m-1}}) = \text{res
Tr}d_{C/K}(\frac{\eta}{z^{m-1}}).$$
And the residue of an exact form is vanishing.
\end{proof}

\begin{lem}$\{c_1(\omega(\sD)), \nabla\}$ is independent of the choice of the
trivializations $s_i$. 
\end{lem}
\begin{proof}First we show independence of the choice of lifting of $s$. As
remarked above,
$g$ is determined by the local lifting of $s$ to $\omega(\sD)$, so
$\{c_1(\omega(\sD)), \nabla\}$ depends only on that choice. If $s$ and $s'$ are
two such local liftings, with $gs=g's'$, we have
$$dgg^{-1}=dg'g'{}^{-1}+d\log(\frac{s'}{s})I= dg'g'{}^{-1}+z^mh
$$
for some $h\in M_r(\sO_C)$. It follows immediately that 
$$\text{res Tr}(dgg^{-1}A) =
\text{res Tr}(dg'g'{}^{-1}A)
$$ 
as desired. 

Next we show independence of the choice of trivializations themselves. Let
$f$ be a rational function on $C$ whose divisor $(f)$ is disjoint from the
singular locus of $\nabla$. It will suffice to show that $s$ and $fs$ as
trivializations give rise to the same invariant, i.e.
\begin{multline}\label{5.5}c_1(\omega(\sD),
s) \cdot {\rm det}( \nabla) - \sum_i{\rm res \
Tr}(dg_i g_i^{-1}A_i) = \\
c_1(\omega(\sD),
fs) \cdot {\rm det}( \nabla) - \sum_i{\rm res \
Tr}((dg_i g_i^{-1}-dff^{-1}I)A_i)
\end{multline}
Recall we can calculate $c_1(\omega(\sD), s) \cdot {\rm det}( \nabla)$ by
choosing $\delta$ a divisor in the linear series $\omega(\sD)$ compatible with
the rigidification $s$ and then restricting $\nabla|_\delta$ and taking the
norm to $\Spec(K)$. Associated to the trivialization $fs$ we may take the
divisor $\delta+(f)$. It follows that
$$c_1(\omega(\sD),
fs) \cdot {\rm det}( \nabla) - c_1(\omega(\sD),
s) \cdot {\rm det}( \nabla) = \text{Norm}\det\nabla|_{(f)}
$$
(To get this relation, one could have taken the formula
\eqref{bigform} as well). 
On the other hand, 
since the formula depends only on the local behavior of $f$
near $\sD$, by
suitably choosing $f$, we may assume $\nabla$ is defined by
$A$ in a neighborhood of $(f)$ and that $f\equiv 1$ modulo some large power of
the maximal ideal at the finite set of points where the connection is not
given by $A$. We can interpret $\text{Tr}
(dff^{-1}A)\in  \OK \otimes \omega_{k(C)}$, so
the sum of the residues over all closed points of $C$ will vanish in $\OK$.
Thus
$$\sum_i\text{res Tr}(dff^{-1}A_i) = -\sum_{(f)}\text{res Tr}(dff^{-1}A) =
\text{Tr}(A|_{(f)})
$$ 
Since the connection matrix for the determinant connection is the trace of the
connection matrix, the contributions to \eqref{5.5} cancel.

\end{proof}

Now consider what happens to the expression
\begin{equation}\{c_1(\omega(\sD)), \nabla\}
\end{equation}
under a change of coordinates given by a matrix $h$ of functions. We have 
\begin{equation}A \mapsto A' := hAh^{-1}+dhh^{-1}=
g'\frac{dz}{z^m}+\frac{\eta'}{z^{m-1}}
\end{equation}
Note that $h$ is regular, so $dhh^{-1}$ does not contribute to the polar part of the
connection, i.e.
\begin{equation}g'=hgh^{-1}+z^ma;\quad \eta' = h\eta h^{-1}+z^{m-1}b
\end{equation}
with $a$ and $b$ regular. 

We compute
\begin{multline}dg'g'{}^{-1}= d(hgh^{-1})hg^{-1}h^{-1}+ez^m+fz^{m-1}dz = \\
 = dhh^{-1}+hdgg^{-1}h^{-1}-hgh^{-1}dhg^{-1}h^{-1}+ez^m+fz^{m-1}dz
\end{multline}
with $e$ and $f$ regular. Thus
\begin{multline}\label{5.11} \text{res Tr}(dg'g'{}^{-1}A') = \text{res
Tr}(dg'g'{}^{-1}hAh^{-1}) = \\ \text{res Tr}(h^{-1}dh A) + 
\text{res Tr}(dgg^{-1}A)- \text{res
Tr}(h^{-1}dhg^{-1}Ag).
\end{multline}
Note
\begin{equation}A-g^{-1}Ag = z^{1-m}(\eta - g^{-1}\eta g) =
z^{1-m}g^{-1}(g\eta - \eta g)
\end{equation}
{From} lemma \ref{lem5.2},2, this expression is regular. Plugging into
\eqref{5.11}, we conclude.
\begin{lem}$\{c_1(\omega(\sD)),\nabla\}$ is independent of the choice of basis
for the bundle $E$.
\end{lem}

\begin{remarks}\begin{enumerate}\item The definition of
$\{c_1(\omega(\sD)),\nabla\}$ was inspired by the calculations in section 4
(cf. formula \eqref{bigform}) for a rank one connection. The formula
\begin{gather}
{\rm det}H^*_{DR}(U, (\sE, \nabla))^{-1} = 
\{c_1(\omega(\sD)), \nabla\}
\end{gather}
when applied to the rank 1 case, gives back the
main theorem of this article. 
\item When $m_i=1$ for all $c_i$, that
is when $\nabla$ has regular singularities, then 
the argument from theorem \ref{thm4.5} (slighlty modified in the
higher rank case) gives
\begin{gather}
\{c_1(\omega(D)), \nabla\} =
c_1(\omega(D), {\rm res}_D) \cdot {\rm det}(\nabla)
\end{gather}
where,  as in theorem \ref{thm4.5}, ${\rm res}_D$ refers to the
natural trivialization coming from the residue. 
\item Finally, twisting $\nabla$ by $f^*\alpha$, where $\alpha \in
\Omega^1_K$ comes from the base, changes the right hand side of
the formula by 
$$(2g-2 + n -\sum_i m_i)r\alpha = -\chi(H_{DR}^*(U,
(\sE, \nabla_{/K})))\alpha ,$$
as it should. Here $r=\rank E$.
\end{enumerate}
\end{remarks}

\begin{prop} \label{prop5.6}
Let $\sD\subset C$ be a divisor on a smooth curve, and let $f\in \sO_{C,\sD}$
be a local defining equation for $\sD$. Let 
$$\nabla : E \to E\otimes\omega(\sD)
$$
be a connection, and write
$$g=\nabla|_\sD : E/E(-\sD) \to (E(\sD)/E)\otimes\omega
$$
Let $j:E-\sD \inj E$ be the inclusion and consider the connection 
$$j_*j^*\nabla : j_*j^*E \to j_*j^*E\otimes \omega.
$$
There is a natural map
$$\iota:\{E\to E\otimes\omega\} \to \{j_*j^*\nabla : j_*j^*E \to
j_*j^*E\otimes \omega\}. 
$$ 
The map $\iota$ is a quasi-isomorphism if and only if for any
$n\ge 0$ the natural map 
$$g-n\cdot id\otimes\frac{df}{f}:E(n\sD)/E((n-1)\sD) \to
(E((n+1)\sD)/E(n\sD))\otimes\omega 
$$
given by $f^{-n}e\mapsto f^{-n}g(e)-nf^{-n-1}e\otimes df$ is a
quasi-isomorphism.

In particular, if,
every point of $\sD$ has multiplicity $\ge 2$, then $g$ is an
isomorphism if and only if $\iota$ is a quasi-isomorphism.
\end{prop}
\begin{proof} The usual exact sequence reduces us to showing the
condition is equivalent to 
$$j_*j^*E/E \stackrel{\cong}{\to}(j_*j^*E/E)\otimes \omega(\sD).  
$$
Writing as usual $\sD=\sD'+D$ where $D$ is the reduced divisor
with support equal to the 
support of $\sD$, we have a commutative square
$$\begin{CD} E/E(-\sD) @>g -n id\otimes \frac{df}{f} >> (E(\sD')/E(-D))\otimes
\omega(D) \\ @V\cong V``\cdot f^{-n}" V @V\cong V ``\cdot f^{-n}" V \\
E(n\sD)/E((n-1)\sD) @>g>>
(E(n\sD+\sD')/E((n-1)\sD+\sD'))\! \otimes\! \omega(D)
\end{CD}
$$
The map $\iota$ is a quasi-isomorphism if and only if for all
$n\ge 0$ the map $g$ on the 
bottom line of the above square is an isomorphism, and this will
hold only if the top line 
is. 

In particular, since $\frac{df}{f}$ has poles of order 1, 
if all multiplicities are
$\geq 2$, then $\iota$ is a quasiisomorphism if and only if $g$
is an isomorphism.
\end{proof}

We close those remarks by a numerical computation 
for $E= \oplus_1^r \sO$ on
$U=\A^1_K$, with parameter $t$ . 
There is only one singular point at $\infty$.
Let us write
$$A = \sum_{i=0}^{m-1}B_it^i + \sum_{i=0}^{m-2} C_i t^idt$$
where the $B_i$ and the $C_i$ 
are matrices with coefficients in $\Omega^1_K$
respectively $K$. The assumption \ref{ass5.3} means that
$C_{m-2} \in GL(r, K)$. We consider the cases $m=2,3$. For $m=2$ both sides
of the  formula are 0, and for $m=3$ they are equal to
\begin{gather}
 {\rm Tr}(B_0 -B_1 C_0C_1^{-1} + B_2 C_0C_1^{-1}C_0C_1^{-1})
\end{gather}
Note that $-C_0C_1^{-1}$ is the zero of the ``polynomial''
$C_0 + C_1t$, where $t$ is a ``variable'' of matrices, and thus
the formula could also be written as 
\begin{gather}
{\rm Tr} A|_{{\rm zero}(C_0 + C_1t)}
\end{gather}
if it made sense.
For higher $m$, the right hand side should be a sort of restriction of
$\nabla$ to the ``Newton'' classes of $\sum_{i=0}^{m-2} C_i t^i$.
 
In the remaining part of this section, we show that the product 
$\{c_1(\omega(\sD)), \nabla\}$ is a particular case of a more
general product between higher rank connections and a larger
class of trivializations along $\sD$.

We consider the tuples $\{E, L, \nabla, \sD, g\}$, where
\begin{gather}
\nabla: E \to E \otimes \Omega^1_C<D>(\sD')
\end{gather}
is a connection on a rank $r$ vector bundle $E$. We denote
by $\nabla|_{\sD}$ the $\sO_C$-linear map
\begin{gather}
\nabla: E \to E \otimes \Big(\Omega^1_C<D>(\sD')/\Omega^1_C\Big).
\end{gather}
Further, $L$ is a rank 1
bundle, $g$ is a trivialization
\begin{gather}
g: E|_{\sD} \cong E\otimes L|_{\sD},
\end{gather}
which fulfills
\begin{gather} \label{comm}
0=[g, \nabla|_{\sD}]: E|_{\sD} \to E\otimes L|_{\sD}
\end{gather}

By lemma \ref{lem5.2}, if $L = \omega_C(\sD)$ and $g$ is the
trivialization  $E|_{\sD} \to E\otimes \Big(\omega_C(\sD)/\omega_C\Big)$ 
arising from the principal part of $\nabla_{C/K}$, then the 
condition \ref{comm} is fulfilled. 

Let us introduce the cocyles of the tuple $\{E,L, \nabla, \sD,
g\}$, as in section 4. If $E$ has cocycle $c_{ij} \in GL(r,
\sO_C)$, $L$ has cocycle $\lambda_{ij} \in \sO^*_C$, $g$ has
cocyle $\mu_i \in GL(r, \sO_{\sD})$, and $\nabla$ has cocyle
$\omega_i \in M(r \times r, \Omega^1_C<D>(\sD'))$ then one has
\begin{gather} \label{coc}
dc_{ij}c_{ij}^{-1} =
\omega_i-c_{ij}\omega_jc_{ij}^{-1} \\
\mu_i=c_{ij}\mu_jc_{ij}^{-1}\lambda_{ij}.
\end{gather}
The commutativity condition \ref{comm} then reads
\begin{gather}\label{comm2}
[\mu_i, \omega_i] = 0.
\end{gather}
Let $\tilde{\mu}_i \in GL(r, \sO_C)$ be a local lifting ot
$\mu_i\in GL(r, \sO_{\sD})$. Then one has
\begin{thm}\label{prod}
The cochain
\begin{equation} \label{coch}
\Big\lbrace\{\lambda_{ij},\det(c_{jk})\},
d\log\lambda_{ij}\wedge \trace(\omega_j),
\trace(-d\tilde\mu_i\tilde\mu_i^{-1}\wedge 
\omega_i)\Big\rbrace 
\end{equation}
is a cocyle in 
\begin{equation}\sK_2 \to \oxxd\to \Big(\oxxd/\oxx\Big)
\end{equation}
and defines a cohomology class $\{c_1(L), g, \nabla\}$ in
\begin{equation}\H^2\Big(C,\sK_2 \to \oxxd\to \Big(\oxxd/\oxx\Big)\Big)
\end{equation}
\end{thm}
\begin{proof}
First
\begin{multline}\delta(d\log(\lambda)\wedge \trace(\omega))(ijk)= \\
d\log(\lambda_{ij})\wedge \trace(\omega_j)+
d\log(\lambda_{jk})\wedge \trace(\omega_k) - 
d\log(\lambda_{ik})\wedge \trace(\omega_k) = \\
d\log(\lambda_{ik})\wedge \trace(\omega_j-\omega_k) - d\log(\lambda_{jk})\wedge
\trace(-\omega_k+\omega_j) = \\
d\log(\lambda_{ij})\wedge \trace(\omega_j-\omega_k) = d\log(\lambda_{ij})\wedge
d\log(\det(c_{jk})).
\end{multline}
Next computing mod the ideal of $\sD$ and so ignoring the tilde, 
\begin{multline} \delta(\trace(-d\tilde\mu_i\tilde\mu_i^{-1}\wedge
\omega_i))(ij)  = \trace(-d\tilde\mu_j\tilde\mu_j^{-1}\wedge
\omega_j + d\tilde\mu_i\tilde\mu_i^{-1}\wedge
\omega_i) = \\
\trace\Big((d\lambda_{ij}\lambda_{ij}^{-1}\cdot I +c_{ij}^{-1}dc_{ij}
-c_{ij}^{-1}d\mu_i\mu_i^{-1}c_{ij}-c_{ij}^{-1}\mu_idc_{ij}c_{ij}^{-1}
\mu_i^{-1}c_{ij}) \wedge\omega_j+ \\
+d\mu_i\mu_i^{-1}\wedge\omega_i\Big) = 
\trace\Big((d\lambda_{ij}\lambda_{ij}^{-1}\cdot I + c_{ij}^{-1}dc_{ij}
-c_{ij}^{-1}d\mu_i\mu_i^{-1}c_{ij}\\ -c_{ij}^{-1}\mu_idc_{ij}c_{ij}^{-1}
\mu_i^{-1}c_{ij}) \wedge(c_{ij}^{-1}\omega_ic_{ij}-c_{ij}^{-1}dc_{ij})+
d\mu_i\mu_i^{-1}\wedge\omega_i\Big)\\
\end{multline}
By our commutation assumption \ref{comm}
\begin{equation}\trace(c_{ij}^{-1}dc_{ij}\wedge
c_{ij}^{-1}\omega_ic_{ij}) =
\trace(c_{ij}^{-1}\mu_idc_{ij}c_{ij}^{-1} 
\mu_i^{-1}c_{ij} \wedge c_{ij}^{-1}\omega_ic_{ij})
\end{equation}
Also terms with no poles (i.e. terms not involving $\omega$) can
be ignored. We get 
\begin{equation} \delta(\trace(-d\tilde\mu_i\tilde\mu_i^{-1}\wedge
\omega_i))(ij)  = d\lambda_{ij}\lambda_{ij}^{-1}\wedge \trace(\omega_i)=
d\lambda_{ij}\lambda_{ij}^{-1}\wedge {\rm Tr}(\omega_j)
\end{equation}
which is what we want. Note the right hand equality holds
because we are computing 
mod forms regular along $\sD$, and ${\rm Tr}(\omega_i)\equiv
{\rm Tr}(\omega_j)$ mod regular 
forms by the cocycle condition \ref{coc}.

It remains to show that our  cocyle \ref{coch}
does not depend on the choice
of the cocycles $\{c, \omega,\lambda, \mu\}$.

If $\lambda_{ij}$ is replaced by $\lambda_{ij}\delta(\nu)_{ij}$,
then the new \ref{coch} differs from the old one by the cochain
$$\delta\Big( \{\nu_i, {\rm det}(c_{ij})\}, d\log \nu_i \wedge
{\rm Tr}(\omega_i)\Big).$$
If $c_{ij}$ is replaced by $\gamma_ic_{ij}\gamma_j^{-1}$, then
$\omega_i$ is replaced by $d\gamma_i \gamma_i^{-1} +
\gamma_i \omega_i \gamma_i^{-1}$, and 
$\mu_i$ is replaced by $\gamma_i \mu_i \gamma_i^{-1}$. The
commutativity relation implies
$d\gamma_i \mu_i \gamma_i^{-1} -\gamma_i \mu_i \gamma_i^{-1}
d\gamma_i \gamma_i^{-1}=0$. Then the new \ref{coch}
differs from the old one by 
\begin{gather}
\{\{\lambda_{ij}, {\rm det}(\delta(\gamma))_{jk}\},
d\lambda_{ij} \wedge \trace d \log (\gamma_i), \notag \\
\rt \Big(\gamma_i^{-1}d\gamma_i \omega_i
  d \mu_i \mu_i^{-1} \gamma_i^{-1}
d\gamma_i  -\mu_i \gamma_i^{-1} d\gamma_i \mu_i^{-1}
\gamma_i^{-1} d\gamma_i -\mu_i\gamma_i^{-1}d\gamma_i
\mu_i^{-1}\omega_i \Big) \}
\end{gather}
 Using the commutativity relation \ref{comm} as well, we see that
this is this expression is the cochain
\begin{gather}
\delta\Big( \{\lambda_{ij}, {\rm det} \gamma_j\}, 0, 0\Big).
\end{gather}
\end{proof}

Now we consider the image under the map $f:C \to \Spec(K)$:
\begin{equation}
f_*( \{c_1(L), g, \nabla\}) \in
\Omega^1_K/d\log K^*
\end{equation}
of $\{c_1(L), g, \nabla\}$. We want to study closedness for forms in
the image of this map.
\begin{lem}\label{lem5.7} Let $\{a,b,c\}$ be a cocycle as in \eqref{coch}
  representing a class in $\H^2(C,\sK_2 \to \Omega^2_C)$, 
with $db=0$. Then
  $df_*\{a,b,c\}= {\rm res}_\sD(dc)\in \Omega^2_K$.
\end{lem}
\begin{proof}
Another representative of $\{a,b,c\}$ in
$$\H^1(C, \Omega^2_C<D>(\sD') \to
\Omega^2_C<D>(\sD')/\Omega^2_C)/d\log (H^1(C, \sK_2))$$
is of the shape $\{0, b + d\log \alpha, c\}$, thus its
derivative in
$$\H^1(C, \Omega^3_C<D>(\sD') \to
\Omega^3_C<D>(\sD')/\Omega^3_C)$$
is of the shape
$\{0, 0, dc\}$ since $db=0$. Then one applies the commutativity
of the diagramm\begin{tiny}
\begin{equation}\begin{CD}\H^1(C,\Omega^2_C\langle D\rangle(\sD') \to
    \Omega^2_C\langle D\rangle(\sD')/\Omega^2_C) @> d >>
    \H^1(C, \Omega^3_C\langle D\rangle(\sD') \to
    \Omega^3_C\langle D\rangle(\sD')/\Omega^3_C) \\
@VVV @VVV \\
\Omega^1_K @>d>> \Omega^2_K 
\end{CD}
\end{equation}
\end{tiny}

\end{proof}
\begin{thm}\label{flat}
If $\nabla^2=0$, that is if $\nabla$ is flat, then 
$f_*( \{c_1(L), g, \nabla\})\in \Omega^1_{K, {\rm clsd}}/d\log
K^*$, that is the image is flat as well. In particular, this is true for
$\{c_1(\omega(\sD)), \nabla\}$, as predicted by conjecture
\ref{conj}.
\end{thm} 
\begin{proof}
Since, under the integrabiltiy assumption, one has in particular
$d(d\log \lambda_{ij} \wedge \trace \omega_j)= 0$, one can apply 
lemma \ref{lem5.7}. 
One has to compute
\begin{gather}
\gamma= -\rt\ d(d\mu \mu^{-1} \omega) \notag \\
= -\rt (d\mu \mu^{-1})^2 \omega + \rt\ d \mu \mu^{-1}d\omega
\end{gather}
We omit the indices since we compute only with one index. Note in the
calculations which follow $\mu$ is regular and $\omega$ has poles
along $\sD$. We write $a\equiv b$ to indicate that the polar parts of
$a$ and $b$ coincide. 

The condition \ref{comm} implies
\begin{gather}
d\mu \omega + \mu d\omega\equiv d\omega \mu - \omega d\mu
\end{gather}
thus
\begin{gather}
d\mu \mu^{-1} \omega + \mu d\omega \mu^{-1}\equiv d\omega - \omega
d\mu \mu^{-1}
\end{gather}
which implies
\begin{gather}
(d\mu \mu^{-1})^2 \omega + d\mu d\omega \mu^{-1}\equiv d\mu
\mu^{-1}d\omega - d\mu \mu^{-1} \omega d\mu \mu^{-1}
\end{gather}
So taking the trace, one obtains
\begin{gather}
2 \trace ( (d\mu \mu^{-1})^2 \omega) \equiv \trace ((d\mu \mu^{-1} -
\mu^{-1}d\mu)d\omega) 
\end{gather}
Now, under the integrability assumption 
\begin{gather} \label{int}
d\omega= \omega \omega
\end{gather}
(using \eqref{comm2}) one obtains
\begin{gather}
2 \trace ((d\mu \mu^{-1})^2 \omega) =0
\end{gather}
Now we consider the other term
\begin{gather}
\gamma'=\trace (d\mu \mu^{-1} (\omega)^2)
\end{gather}
Choosing a local basis of the bundle $E$, we write $\mu$ as a matrix
$M$, and $\omega$ as 
matrix $\Omega$.
Then the condition \ref{comm}
reads
\begin{gather} \label{com2}
M \Omega= \Omega M
\end{gather}
and one has
\begin{gather}
\gamma'= \trace  (dMM^{-1}\Omega^2)
\end{gather}
The condition \ref{com2} implies
\begin{gather}
\Omega^2 = M^{-1} \Omega^2M
\end{gather}
thus
\begin{gather}
\gamma'= \trace ( dM \Omega^2 M^{-1} )= \trace (M^{-1}dM \Omega^2)
\end{gather}
On the other hand, differentiating the condition \ref{com2}, one
obtains
\begin{gather}
M^{-1} dM \Omega^2 + d\Omega\Omega =\notag\\
M^{-1} d\Omega M \Omega - M^{-1} \Omega dM \Omega= \notag \\
M^{-1}d\Omega \Omega M -\Omega M^{-1}dM\Omega 
\end{gather}
So taking the trace, one obtains
\begin{gather}
\gamma'=\trace  (M^{-1} dM \Omega^2) = 
-\trace  (M^{-1}dM \Omega^2 )= -\gamma'
\end{gather}
Thus $\gamma'=0$.

\end{proof}

\section{Rank $1$ Irregular Connections in Arbitrary Dimension
on Projective Manifolds}

Let $X$ be a smooth, projective variety in characteristic 0, 
and let $D\inj X$ be a normal
crossings divisor. Given $m_i\ge 0$, define $C(X,D,\underline{m})$ to be
the group of isomorphism classes of line bundles $\sL$ on $X$ together with
an integrable connection
$$\nabla : \sL \to \sL\otimes \Omega^1_X<D>(\sD)
$$
where $\sD= \sum_i m_iD_i$.
Define for $m_i\ge 1$
$$\text{Irreg}(X,D,\underline{m}) := C(X,D,\underline{m})/C(X,D,\underline{0}).
$$
\begin{thm}\label{thm6.1} With notation as above, there is a canonical
isomorphism
$$\text{Irreg}(X,D,\underline{m}) \cong \Gamma(X,\sO_X(\sD)/\sO_X),
$$
where $\sO_X(\sD)/\sO_X$ 
is the sheaf of principal parts of degree $\le m_i$ along
$D_i$. 
\end{thm}
\begin{lem}Exterior differentiation induces an isomorphism
$$d: \sO_X(\sD)/\sO_X \stackrel{\cong}{\to} \Omega^1_X<
D>(\sD)_{{\rm clsd}}/
\Omega^1_X< D>_{{\rm clsd}}
$$
\end{lem}
\begin{proof}[proof of lemma] We first check injectivity. Let $x=0$ be a local
defining equation for $D$. Suppose for some $n$ with $1\le n\le m$ we had
$$d(\frac{a}{x^n}) = \frac{1}{x^n}(da - na\frac{dx}{x})
 \in \Omega^1_X<D>
$$
Multiplying by $x^n$ and taking residue along $D$, it would follow that
$a|_D=0$, i.e. $\frac{a}{x^n}=\frac{b}{x^{n-1}}$. 

To show surjectivity, write a local section of 
$\Omega^1<D>(nD)_{{\rm clsd}}$
(here $1\le n\le m$) in the form
\begin{equation}\label{1} \omega = \frac{adx}{x^{n+1}}+\frac{B}{x^n}
\end{equation}
where $B$ does not involve $dx$. Replacing $\omega$ by
$\omega+d(\frac{a}{nx^n})$, we can assume
$$\omega = \frac{adx+B}{x^n}.
$$
Then
$$0 = d\omega = \frac{da\wedge dx + dB - n\frac{dx}{x}\wedge B}{x^n}.
$$
Multiplying by $x^n$ and taking residue along $D$, we see $B|_D = 0$. Since $B$
does not involve $dx$, it follows that $B=xC$ and $\omega$ can be written
$$\omega = \frac{adx}{x^n}+\frac{C}{x^{n-1}}
$$
Comparing with \eqref{1}, we have lowered the order of pole by $1$. This
process continues until $\omega$ has log poles. 
\end{proof}

\begin{proof}[proof of theorem \ref{thm6.1}] Using the lemma, we get a diagram
with exact rows
\begin{equation*}\minCDarrowwidth.5cm \begin{CD} 0 @>>> \sO_X^* @= \sO_X^*
@>>> 0 @>>> 0\\ 
@. @VVV @VVV @VVV \\
0 @>>> \Omega^1_X<D>_{{\rm clsd}} @>>> 
\Omega^1_X<D>(\sD)_{{\rm clsd}} @>>>
\sO_X(\sD)/\sO_X @>>> 0
\end{CD} 
\end{equation*}
We view this as a diagram of complexes written vertically. Using the standard
hypercohomological interpretation of line bundles with connection, this yields
an exact sequence
\begin{multline*}0 \to C(X,D,0) \to C(X,D,\underline{m}) 
\to H^0(X,\sO_X(\sD)/\sO_X) \\
\stackrel{\delta}{\to} \H^2(X,\sO_X^* \to \Omega^1_X<D>_{{\rm clsd}}) 
\end{multline*}
We claim the map $\delta$ above is zero. In the derived category, $\delta$
factors
\begin{multline*}\sO_X(\sD)/\sO_X \stackrel{\cong}{\to} \frac{\Omega^1_X<
D>(\sD)_{{\rm clsd}}}{\Omega^1_X< D>_{{\rm clsd}}} 
\stackrel{\partial}{\to} \Omega^1_X(<D>)_{{\rm clsd}}[1] \\
\to \{\sO^*_X \to \Omega^1_X(<D>)_{{\rm clsd}} \}[2].
\end{multline*}
We have a factorization of $\partial$:
\begin{multline*}H^0(\sO_X(\sD)/\sO_X) \to H^1(\sO_X) \stackrel{\tilde d}{\to}
H^1(\Omega^1_{X, {\rm clsd}}) \to H^1(\Omega^1_X<D>_{{\rm clsd}}),
\end{multline*}
so it suffices to show the map $\tilde d$ is zero. By Hodge theory, the
composition
$$H^1(\sO_X) \stackrel{\tilde d}{\to} H^1(\Omega^1_{X, {\rm clsd}})
\stackrel{\iota}{\to} \H^1(\Omega^1_X \to \Omega^2_X)
$$
is zero, and the map $\iota$ is injective
as the complex $\{\Omega^1_X/\Omega^1_{X,{\rm clsd}} \to
\Omega^2_X\}$ is quasi-isomorphic to the complex
$\{0 \to \Omega^2_X/\Omega^2_{X, {\rm exact}}\}$, and in particular
starts in degree 1.
\end{proof}

\newpage
\bibliographystyle{plain}
\renewcommand\refname{References}

\end{document}